\documentclass[twoside, 11pt]{article}

\usepackage{amssymb, amsmath, latexsym, mathrsfs, verbatim, calc}
\usepackage{cite}
\usepackage{amssymb}
\usepackage[usenames,dvipsnames]{color}
\usepackage{enumerate}
\usepackage{graphicx}
\usepackage[colorlinks, linkcolor=blue, anchorcolor=blue, citecolor=blue]{hyperref}

\numberwithin{equation}{section}
\newtheorem{thm}{Theorem}[section]
\newtheorem{lem}[thm]{Lemma}
\newtheorem{cor}[thm]{Corollary}
\newtheorem{prop}[thm]{Proposition}
\newtheorem{rem}[thm]{Remark}

\newtheorem{defn}[thm]{Definition}
\newtheorem{assum}[thm]{Assumption}

\definecolor{blue-violet}{rgb}{0.54, 0.17, 0.89}
\definecolor{purple1}{rgb}{0.63, 0.36, 0.94}
\definecolor{purple2}{rgb}{0.87, 0.0, 1.0}
\definecolor{purple}{rgb}{0.5, 0.0, 0.5}
\definecolor{pansypurple}{rgb}{0.47, 0.09, 0.29}
\definecolor{orange}{rgb}{1.0, 0.27, 0.0}

\def\ba{\begin{array}}
\def\ea{\end{array}}
\def\beq{\begin{equation}}
\def\endeq{\end{equation}}
\def\bes{\begin{equation*}}
\def\ees{\end{equation*}}
\def\bea{\begin{eqnarray}}
\def\eea{\end{eqnarray}}
\def\beaa{\begin{eqnarray*}}
\def\eeaa{\end{eqnarray*}}

\def\yw{\tilde}
\def\dis{\displaystyle}

\def\no{\noindent}

\def\lastline{\par \vspace{-7.3ex} \no}

\def\nts{\negthinspace}

\def\q{\quad}
\def\qq{\qquad}
\def\pw{\partial}
\def\dy{\triangleq}
\def\gh{\sqrt}
\def\ol{\overline}
\def\ul{\underline}

\def\={=\nts \nts=\nts \nts=\nts \nts=}

\def\bh{\mbox}
\def\({\textnormal{(}}
\def\){\textnormal{)}}

\def\qx{\wedge}

\def\a{\alpha}

\def\d{\delta}
\def\e{\varepsilon}

\def\k{\kappa}
\def\l{\lambda}

\def\r{\rho}
\def\si{\sigma}
\def\t{\tau}
\def\f{\varphi}
\def\th{\theta}
\def\o{\omega}

\def\f{\phi}
\def\vf{\varphi}
\def\p{\psi}
\def\ksfl{\begin{cases}}
\def\jsfl{\end{cases}}

\def\G{\Gamma}
\def\L{\Lambda}
\def\O{\Omega}

\def\P{\Psi}

\def\kspl{\begin{itemize}}
\def\jspl{\end{itemize}}

\def\wq{\infty}
\def\qd{\backslash}

\def\bbh{\nonumber}

\def\bbh{\nonumber}

\def\no{\noindent}

\def\q{\quad}
\def\qq{\qquad}

\def\jf{\int}
\def\fs{\frac}
\def\yjt{\rightarrow}
\def\ysy{\subset}
\def\yj{\hat}
\def\bdy{\neq}
\def\dqy{\hx{D}}
\def\kspl{\begin{enumerate}[(1)]}
\def\jspl{\end{enumerate}}
\def\zb{\left}
\def\yb{\right}
\def\hx{\mathscr}
\def\qed{\hfill \rule[0cm]{.25cm}{.25cm}\medskip}   

\def\b1{{\bf 1}}

\newenvironment{itm}{\vspace{-1ex}\begin{itemize}}{\end{itemize}}
\def\bi{\begin{itm}}
\def\ei{\end{itm}}

\def\equ_ind{\arabic{section}.\arabic{equation}}
\def\sec_ind{\arabic{section}}

\usepackage{CJK}
\begin{document}
\begin{CJK*}{GBK}{song}
\title{\Large \bf Optimal Singular Dividend Problem under the Sparre Anderson Model
\thanks{This research is supported by Chinese NSF grants No. 11471171 and No. 11571189.}}

\author{{Linlin Tian}$^a$\footnote{E-mail:linlin.tian@mail.nankai.edu.cn}~~~
{Lihua Bai}$^a$\footnote{Corresponding author. Phone number: (+86)13920931661. E-mail:lhbai@nankai.edu.cn }~~~ {Junyi Guo}$^a$\footnote{E-mail:jyguo@nankai.edu.cn}
\\
 \small  a. School of Mathematical Sciences, Nankai University, Tianjin 300071, China.}\,
\date{}
\maketitle
\begin{center}
\begin{minipage}{130mm}
{\bf Abstract.}
Consider an insurance company for which the reserve process follows the Sparre Anderson model. In this paper, we study the optimal dividend problem for such  a company as Bai, Ma and Xing \cite{bai2017optimal} do. However, we remove the constant restriction on the dividend rates, i.e. the optimization problem is of singular type. In this case, the value function is no longer bounded and the associated HJB equation is a variational inequality involving a first order integro-differential operator and a gradient constraint. We  use other techniques to prove the regularity properties for the value function and
show that the value function is a constrained viscosity solution of the associated HJB equation. In addition, we show that the value function is the upper
semi-continuous envelop of the supremum for a class of subsolutions.

\vspace{3mm} {\bf Keywords:}  Singular control, optimal dividend, HJB equation, constrained viscosity solution, viscosity supersolution, viscosity subsolution.

\vspace{3mm} {\bf 2010 Mathematics Subject Classification}: 34H05, 34K35, 93E20, 49L20, 49L25.
\end{minipage}
\end{center}

\section{Introduction}
The dividend problem was first posed by De
Finetti  \cite{de1957impostazione} at the International Congress of Actuaries in 1957. Asmussen and Taksar\cite{Asmussen1997Controlled}  solved the optimal dividend problem for the special case of Brownian motion. They determined that the optimal strategy is a constant barrier strategy in the case of  unbounded dividend and a so-called threshold strategy in the case of restricted
dividend rates. 
In the case of a surplus process following a compound Poisson process,
Gerber and Shiu \cite{gerber2006optimal} showed that the optimal strategy is a threshold strategy if claim sizes are
exponentially distributed for restricted dividend rates. Later, Azcue and Muler\cite{azcue2005optimal} considered the problem of maximizing the cumulative expected discounted dividend payouts of the insurance company where the reserve process follows the Cram$\acute{\bh{e}}$r-Lundberg model. In the setting of a jump-diffusion process,  Belhaj\cite{belhaj2010optimal} determined that  the optimal dividend policy is a barrier strategy if claim jumps are exponentially distributed.     Optimal dividend control problems have been extensively researched, such as in  Albrecher and Thonhauser\cite{albrecher2009optimality} and the exhaustive collection of references cited therein   for the past  development of research into such problems.

It is well known that the closed form of the value function for the optimal dividend problem is difficult to obtain if reserve follows the compound Poisson process with a general claim distribution. In this case, Azcue and Muler \cite{azcue2005optimal} investigated the optimal dividend-reinsurance problem and then they studied the optimal dividend-investment problem in \cite{azcue2010optimal}. Both papers used the notion of a viscosity solution to construct the connection between   the value function and the HJB equation. If the compound Poisson claim process  is replaced by  a renewal process, the result is a model known as Sparre Andersen risk model \cite{andersen1957collective};   we should also study the dividend problem in the framework of  viscosity solution. As the dividend problem under the renewal process is non-Markovian, it is more challenging to investigate.  Li and Garrido \cite{li2004class} computed an IDE for the Gerber-Shiu function in the case of a barrier strategy; 
Albrecher et al. \cite{albrecher2005distribution} calculated the moments of the expected discounted dividend payments under a barrier strategy;  Albrecher and  Hartinger \cite{albrecher2006non} showed
that even in the case of Erlang(2) distributed interclass times and exponentially distributed
claim amounts, a horizontal barrier strategy is not optimal anymore, as it can be outperformed
by a strategy that depends on the time elapsed since the previous claim occurrence.

We consider the optimal dividend problem of the renewal process in the framework of stochastic optimal control.  There is a rigorous connection between stochastic optimal control and fully nonlinear integro-partial differential equations. For instance, Benth, Karlsen and Reikvam\cite{benth2001optimal} studied the optimal portfolio selection problem using the viscosity solution of an integro-partial differential equation. Seydel \cite{seydel2009existence}
studied the  optimal impulse control problems for the compound Poisson jump diffusion process using the viscosity solution of a quasi-variational inequality. For the renewal process,
Bai, Ma and Xing \cite{bai2017optimal} studied the optimal dividend problem and investment problem under the Sparre Anderson Model with a constant restriction $M<\wq$ on the dividend rates. Based on \cite{bai2017optimal}, we explore the optimal singular dividend problem under the Sparre Andersen model, i.e., the jump dividend is allowed in our model. The difference between this paper and \cite{bai2017optimal} primarily consists of three aspects:
\begin{itemize}
 \item   In the context of showing the continuity of the value function with respect to variable $x$, \cite{bai2017optimal} introduced a penalty function, while we consider this question by constructing admissible strategies. This also leads to a difference in the proof of continuity of the value function with respect to variable $w$.
 \item  As our HJB equation is a first-order integro-differential equation  with a gradient constraint,   while the HJB equation of \cite{bai2017optimal} is a  second-order integro-differential equation,  the proof of verification theorem  is  different from that of \cite{bai2017optimal}.
     \item The structure of the HJB equation and the boundary condition of the value function are changed compared to those of \cite{bai2017optimal}. Based on the HJB equation and the boundary condition, we construct a viscosity subsolution and a supersolution of the HJB equation and we provide a candidate of the value function.
\end{itemize}
The paper is organized as follows.  First, we prove several continuity properties of the value function and develop the dynamic programming principle for the optimal dividend problem. Second, we characterize the value function of the singular dividend problem as a constrained viscosity solution of the associated Hamilton-Jacobi-Bellman equation. Finally, we show that the value function is the upper semi-continuous envelop of the supremum for a class of subsolutions.
\section{Model and Assumptions}
Let $(\Omega,\mathbb{P},\hx{F})$ be a probability space and $(\hx{F}_t)$ be a given filtration satisfying the usual assumptions. We consider a renewal counting process $N=\{N_t\}_{t\ge 0}$ on $(\Omega,\mathbb{P},\hx{F})$. Denote $\{\si_n\}_{n\ge 0}$ as the jump times of the renewal  counting process $N$, and denote $T_i=\si_i-\si_{i-1},i=1,2,\cdots$ as the waiting times between claims. We also assume that $\{T_i\}_{i=1}^{\wq}$ are  independent and identically distributed with the common distribution $F:\mathbb{R}_+\yjt [0,1]$, and
that there exists an intensity function $\l:[0,\wq)\yjt [0,\wq)$ such that $\bar{F}(t):=\mathbb{P}\{T_1>t\}=\exp\{-\jf_0^t \l(u)du\}$. In other words, $\l (t)=f(t)/\bar{F}(t), t\ge 0$, where $f$ is the common density function of $T_i$'s.
An important feature of the Sparre Anderson model is the ``compound renewal process''  that describes the total claim process $\sum_{i=1}^{N_t}U_i$, where $N_t$ is the renewal process representing the frequency of claims up to time $t$, and $\{U_i\}_{i=1}^{\wq}$ is a sequence of random variables representing the size of incoming claims. We assume that $\{U_i\}$ are independent, identically distributed with the common distribution $G:\mathbb{R}^+\yjt  [0,1]$, and are independent of $N$. Denote $Q_t=\sum_{i=1}^{N_t}U_i, t\ge 0$. As the process $Q_t$ is non-Markovian in general (unless the counting process is a Poisson process), we cannot use the dynamic programming approach directly. Instead, we use the so-called $Backward\; Markovization$ technique that was used in \cite{bai2017optimal}. Thus, we need to introduce a new process
\[W_t=t-\sigma_{N_t},t\ge0,\]
of the time elapsed since the last jump. We observe that $0\le W_t\le t\le T$.  It is known that the process $(t,Q_t,W_t),t\ge 0$ is a piecewise deterministic Markov process (see, e.g.,~\cite{rolski1998stochastic}). Throughout this paper, we consider the following filtration $\{\hx{F}\}_{t\ge 0}$, where $\hx{F}=\hx{F}^N \vee \hx{F}^W$, $t\ge 0$. Here, $\hx{F}^N,\hx{F}^W$ denote the natural filtrations generated by processes $N,W$, respectively, with the usual $\mathbb{P}$-augmentation such that it satisfies the $usual$  $hypotheses$; see, e.g.,~\cite{Protter1990Stochastic}.

An important feature of the dynamic optimal control theory is to allow  the starting point to be any time $s\in[0,T]$. In our Sparre Anderson model, we consider the initial time $s\in [0,T]$ with the initial elapsed time $W_s=w$   instead of $t=0$.  We consider the $regular$ $conditional$ $probability$ $distribution$ $\mathbb{P}_{sw}(\cdot)=\mathbb{P}[\cdot|W_s=w]$ on $(\Omega, \hx{F})$ and consider the ``shifted'' version of process $(Q,W)$ on space $(\Omega, \hx{F}, \mathbb{P}_{sw};\mathbb{F}^s)$, where $\mathbb{F}^s=\{\hx{F}\}_{t\ge s}$.
Next, we restart the clock at time $s\in [0,T]$ by defining the new counting process $N_t^s:=N_t-N_s, t\in [s,T]$.
Then, under $P_{sw}$, $N^s$ is a ``delayed'' renewal process, in the sense that while its waiting times $T_i^s,i\ge 2$   remain independent and identically distributed  with the original processes $T_i$'s, its ``time-to-first jump'',
denoted by $T_1^{s,w}\dy T_{N_s+1}-w=\si_{N_s+1}-s$, should have the survival probability
\[\mathbb{P}_{sw}\{T_1^{s,w}>t\}=\mathbb{P}\{T_1>t+w|T_1>w\}= \exp\left[-\jf_w^{w+t}\l(u)du\right].\]
To emphasize the dependence of $N$ on $w$,   in the following, denote $N_t^s|_{W_s=w}\dy N_t^{s,w}$, $t\ge s.$ We denote $Q_t^{s,w}=\sum_{i=1}^{N_t^{s,w}}U_i$ and $W_t^{s,w}= w+W_t-W_s$, $t\ge s$.  It is easy to see that $(Q_t^{s,w},W_t^{s,w}),t\ge s$, is an $\mathbb{F}^s$-adapted Markov  process defined on $(\O, \hx{F}, \mathbb{P}_{sw})$. We will rely on the following standing assumptions.
\begin{assum}
\begin{enumerate}[(1)]\label{a52}
\item The insurance premium $p$ and discount factor $c$ are all positive constants.
\item The distribution function $F$ (of $T_i$'s) and $G$ (of ${U_i}$'s) are continuous on $[0,+\wq)$. The intensity function of $F$ is denoted by $f(t)$; there exists a continuous and bounded intensity function such that $\l(t)=f(t)/\bar{F}(t)>0, t\in [0,T]$. There exists a constant $\L>0$ such that for all $w\in[0,T]$, it holds that $\l(w)\le\L$.
\end{enumerate}
\end{assum}

From now on, we consider the optimal dividend problem. For a given $s\in [0,T]$ and any dividend policy $L_t^{s,w}$, where $L_t^{s,w}$ denotes the cumulative dividend from time $s$ to time $t$, the SDE of the wealth process  $X_t^{\pi,x}$ satisfies
\[X_t^{\pi,x}=x+p(t-s)-Q_t^{s,w}-L_t^{s,w}.\]
We denote $(X^{\pi,x}, W,L)=(X^{\pi,s,x,w}, W^{s,w},L^{s,w})$ for simplicity.
We call a control strategy $\pi=\{L_t\}_{t\ge 0}$ admissible, if the following hold true:\\
(1) It is $\hx{F}-$predictable, nondecreasing, and  c\`{a}gl\`{a}d.\\
(2) For any time $t\ge 0$, the process $L_t$ satisfies
\bea\label{a48}
L_t\le x+p(t-s)-Q_t^{s,w},
\eea
which means that the dividend process $L_t$ cannot cause  bankruptcy.
Denote  $U_{ad}^{x,w}[s,T]$ as the set of all admissible control strategies with initial wealth $x$ and elapsed time $w$ since the last jump at time $s$. From (\ref{a48}), we observe that the set of admissible strategies is related to initial capital $x$, which is quite different from \cite{bai2017optimal}.
For the given initial date $(s,x,w)$, we define the cost function by
\[J(s,x,w;\pi)\!=\!\mathbb{E}_{sxw}\!\left[\int_{s}^{({\tau^{\pi,x}}\land T)+}\!\!\!e^{-c(t-s)}dL_t^\pi\yb]
\!\triangleq\!\mathbb{E}\!\left[\int_{s}^{({\tau^{\pi,x}}\land T)+}e^{-c(t-s)}dL_t|X_s=x, W_s=w\right],\]
where $\t^{\pi,x}=\inf\{t\ge s: X_t^{\pi,x}<0\}$ denotes the ruin time of the insurance company,   and $c>0$ is the discounting factor. From now on, denote $X_t^\pi= X_t^{\pi,x} $ and $\t^\pi=\t^{\pi,x} $ for simplicity, if such notation does not cause confusion. $J(s,x,w;\pi)$ is the expected total discounted amount of dividends before ruin. Our objective is to find the optimal strategy to maximize the expectation of cumulative discounted dividends. The value function is defined by
\[V(s,x,w)=\sup_{\pi\in U_{ad}^{x,w}[s,T]}J(s,x,w;\pi).\]
The function $V(s,x,w)$ should be defined on $\{(s,x,w):0\le s\le T, x\ge 0, 0\le w\le s\}$. Thus, we introduce several notations   for simplicity. Denote
\begin{align}\nonumber
D&=\{(s,x,w)|0\le s\le T, 0\le x,0\le w\le s\},\\\bbh
\hx{D}&=\bh{int} D=\{(s,x,w)\in D| 0<s<T,0<x, 0<w<s\},\\\bbh
\hx{D}^*&=\{(s,x,w)\in D|0\le s<T,0\le x,0\le w\le s\}.
\end{align}
Here, we note that $\hx{D}^*$ does not include the boundary of $s=T$.

\begin{rem}
As the singular dividend policy is admissible in our model, it may have a jump at any time $t\le T$. Thus, the optimal choice at time $T$ is to pay the entire wealth as dividends.
\end{rem}
\section{Basic Properties of The Value Function}
In this section, we present several propositions to characterize a number of regularity properties of the value function $V(s,x,w)$.
\begin{prop}\label{25}  Under assumption \ref{a52},
the optimal return function $V$ is well-defined on $D$, and for all $x\ge 0, s\in[0,T],w\in [0,s]$, it holds that
$x\le V(s,x,w)\le x+p(T-s)$.
\end{prop}
$\mathbf{Proof}$ \q
From $E\zb[\jf_s^{(\t^\pi\qx T)+}e^{-c(t-s)}dL_t\yb]\le L_{\t^\pi\qx T}$, we observe that for any strategy $\pi\in U_{ad}^{x,w}[s,T]$, it holds that
$
 J(s,x,w;\pi)\le \mathbb{E}_{sxw}\zb[L_{(\t^\pi\qx T)+}\yb].
$
By (\ref{a48}), we observe that $L_{(\t^\pi\qx T)+}\le x+p(T-s)$. On the other hand, the second proposition is trivial. As we allow for the company to have a jump dividend, the company can pay the entire capital wealth $x$ at initial time $s$ as dividends and then pay  dividends at rate $p$ after time $s$. Thus, we obtain $V(s,x,w)\ge x$.  $\hfill\Box$
\begin{prop}\label{a50}
Under assumption \ref{a52}, for all $(s,x,w),(s+h,x,w)\in D$, $h>0$, it   holds that

 \no
(1) $ V(s+h,x,w)-V(s,x,w)\leq 0$.   \\(2) $V(s,x,w)-V(s+h,x,w)\le 2ph $.
\end{prop}
$\mathbf{Proof}$\q (1) For the first conclusion, the proof is identical to that of Proposition  3.3 in \cite{bai2017optimal}. We omit it here.

(2) For any strategy $\pi\in{U}_{ad}^{x,w}[s,T]$, 
we construct $\pi^h\in U_{ad}^{x,w}[s,T]$ as follows:
\[\pi_t^h=1_{\{\tau^{\pi}>T-h\}}\zb[1_{\{t< T-h\}}\pi_t+1_{\{t=T-h\}}\pi_t^{all}+1_{\{t>T-h\}}\pi_t^p\yb]+1_{\{\tau^\pi\le T-h\}}\pi_t,\]
where $\pi^{all}$ means  using  the entire wealth as dividend immediately, and $\pi^p$ means paying dividend at rate $p$. Denote $\t^h$ as the ruin time of strategy $\pi^h$ and $L_t^h$ as the corresponding cumulative dividend payments of $\pi^h$.   Thus,
\begin{eqnarray}
\nonumber &&J(s,x,w;{{\pi}^h})
=\mathbb{E}_{sxw}\left[\int_{s}^{(\t^h\wedge T)+}e^{-c(t-s)}dL_t^h\right]\\\nonumber &=&\mathbb{E}_{sxw}\left[\int_{s}^{({\t^h}
\wedge T)+}e^{-c(t-s)}dL_t^h:\t^h>T-h\right]\!\!+\!\mathbb{E}_{sxw}\left[\int_{s}^{{\t^h}\wedge T}\!\!\!
e^{-c(t-s)}dL_t^h:\t^h\!\leq\! T\!-\!h\right]\\\bbh&=&\mathbb{E}_{sxw}\left[\int_{s}^{(T-h)+}\!\!\!e^{-c(t-s)}dL_t^h
:\t^h\!>\!T\!-\!h\right]\!\!+\!\!\mathbb{E}_{sxw}\left[\int_{(T-h)+}^{({\t^h}\wedge T)+}\!\!\!e^{-c(t-s)}dL_t^h
:\t^h>T-h\right]\\\label{a19}& &+\mathbb{E}_{sxw}\!\left[\int_{s}^{{\t^h}}\!\!\!e^{-c(t-s)}dL_t^h:{\t^h}\leq T-h\right].
\end{eqnarray}
Note that  $\pi^h|_{[s,T-h]}$ can be regarded as an admissible strategy of ${U}_{ad}^{x,w}[s+h,T]$. Thus, we obtain that
\begin{align}\label{a18}
\mathbb{E}_{sxw}\left[\int_{s}^{(T-h)+}\!\!\!\!e^{-c(t\!-\!s)}dL_t^h
:\t^h\!>\!T\!-\!h\right]\!\!+\!\mathbb{E}_{sxw}\left[\int_{s}^{{\t^h}}\!\!e^{-c(t-s)}dL_t^h:{\t^h}\!\leq \!T\!\!-\!\!h\right]\!
\!\le\! V(s\!+\!h,x,w).
\end{align}
Substituting (\ref{a18}) into (\ref{a19}), and considering that on $\{\t^\pi>T-h\}$, $X_{(T-h)+}^h=0$,
and all dividends of $\pi^h$ during the interval $(T-h,T]$ are less than $ph$,  we obtain \begin{align}\label{a23}
J(s,x,w;{{\pi}^h})\!\le\! V(s\!+\!h,x,w)\!\!+\!\mathbb{E}_{sxw}\!\!\left[\int_{(T-h)+}^{{\t^h}\wedge T}\!\!\!e^{-c(t-s)}dL_t^h:\!\t^h\!>\!T\!-\!h\right]\!\!\le\! V(s\!+\!h,x,w)\!+\!ph.\end{align}
In what follows, we consider the difference between $\pi^h$ and $\pi$. By the definition of $\pi^h$,  we can calculate the following:
\begin{align}\bbh
&J(s,x,w;{\pi})-J(s,x,w;{\pi^h})=\mathbb{E}_{sxw}\left[\int_{s}^{({\tau^\pi}\wedge T)+}e^{-c(t-s)}dL_t^\pi-\int_{s}^{({\t^h\wedge  T})+}e^{-c(t-s)}dL_t^h\right]\\\bbh
&=\mathbb{E}_{sxw}\left[\int_{s}^{(\tau^{\pi}\wedge T)+}e^{-c(t-s)}dL_t^\pi-\int_{s}^{(\t^h\wedge T)+}e^{-c(t-s)}dL_t^h:\tau^\pi\leq T-h\right]\\\bbh
&\q+\mathbb{E}_{sxw}\left[\int_{s}^{(\tau^\pi\wedge T)+}e^{-c(t-s)}dL_t^\pi -\int_{s}^{(\t^h\wedge T)+}e^{-c(t-s)}dL_t^h:\tau^\pi>T-h\right].
\end{align}

As there is no difference between $\pi^h$ and $\pi$ on $[s,T\!-\!h)$ on set $\{\t^\pi\!>\!T-h\}\cup\{\t^\pi\le \!T-\!h\}$, we obtain
\begin{align}
\label{a20}
J(s,x,w;{\pi})\!\!-\!\!J(s,x,w;{\pi^h})
\!=\!\mathbb{E}_{sxw}\left[\int_{T-h}^{(\tau^\pi\wedge T)+}\!\!\!e^{-c(t-s)}dL_t^{\pi}\!\!-\!\!\int_{T-h}^{(\t^h\qx T)+}\!\!\!e^{-c(t-s)}dL_t^h:\tau^{\pi}\!\!>\!\!T\!\!-\!\!h\right].
\end{align}
As the dividend strategy $\pi$ cannot cause bankruptcy, we observe that
\begin{align}\label{a21}
\mathbb{E}_{sxw}\left[\int_{(T-h)-}^{(\tau^\pi\wedge T)+}e^{-c(t-s)}dL_t^{\pi}:\tau^{\pi}>T-h\right]\le \mathbb{E}_{sxw}\left[ph+X_{T-h}^{\pi}:\tau^{\pi}>T-h\right].
\end{align}
By the definition of $\pi^h$, the dividend $\pi^h$ will have a jump at time $T-h$ on $\{\tau^{\pi}>T-h\}$; thus,
\begin{align}\label{a22}
\mathbb{E}\left[\int_{(T-h)-}^{(\t^h\qx T)+}e^{-c(t-s)}dL_t^h:\tau^{\pi}>T-h\right]\ge \mathbb{E}_{sxw}\left[ X_{T-h}^\pi :\tau^{\pi}>T-h\right].
\end{align}
Substituting (\ref{a21}) and (\ref{a22}) into (\ref{a20}), we obtain
\begin{align}\label{a24}
J(s,x,w;{\pi})-J(s,x,w;{\pi^h})\le \mathbb{E}_{sxw}\left[ph+X_{T-h}^{\pi}-X_{T-h}^{\pi}:\tau^{\pi}>T-h\right]\le ph.
\end{align}
Combining (\ref{a23}) and (\ref{a24}), we obtain that   for any strategy $\pi\in{U}_{ad}^{x,w}[s,T]$,
\[J(s,x,w;\pi)\!\!-\!\!V(s+h,x,w)\!\!=\!\!J(s,x,w;\pi)\!\!-\!\!J(s,x,w;{\pi^h})\!+\!J(s,x,w;{\pi^h})\!\!-\!\!V(s+h,x,w)\!\leq\!2ph.\]
Taking the supremum over all strategies $\pi\in{U}_{ad}^{x,w}[s,T]$, we obtain
\begin{eqnarray}\label{a28}
V(s,x,w)-V(s+h,x,w)\le 2ph.
\end{eqnarray} $\hfill\Box$
\begin{prop} \label{a38} Under assumption \ref{a52}, the following holds:\\
(1) For any $x_1\ge x_2\ge 0$, $V(s,x_1,w)-V(s,x_2,w)\ge x_1-x_2$;\\
(2) For any compact set $\mathfrak{D}$, the mapping $x\mapsto V(s,x,w)$ is continuous, uniformly for $(s,x,w)\in \mathfrak{D}$.
\end{prop}
$\mathbf{Proof}$\q (1) For the first part, it is trivial to note that for any strategy $\pi_2\in U_{ad}^{x_2,w}[s,T]$, we can construct $\pi_1$ such that $\pi_1$ has a jump  dividend $x_1-x_2$ at time $s$ and then follows the strategy $\pi_2$. Thus, we observe that
\[V(s,x_1,w)\ge J(s,x_1,w;\pi_1)=J(s,x_2,w;\pi_2)+x_1-x_2.\]
As $\pi_2\in U_{ad}^{x_2,w}[s,T]$ is arbitrary, taking the supremum over $U_{ad}^{x_2,w}[s,T]$, we obtain
\[V(s,x_1,w)-V(s,x_2,w)\ge x_1-x_2.\]
(2) Here, we borrow several arguments from \cite{scheer2011optimal}.
Suppose that $h\ge 0$, and let $\pi$ be a strategy with initial capital $x$. Now, we construct another strategy $\yw{\pi}$ that pays no dividend until some stopping time $\t^h$ to be defined below. At time $\t^h$, the strategy $\yw{\pi}\in U_{ad}^{x-h,w}[s,T]$ is adjusted so that two surplus processes $X_t^{\pi,x}$ and $X_t^{\yw{\pi},x-h}$ coincide from $\t^h$ on. Let
\[ \t^h\dy\inf\left\{t\ge s: L_t^\pi\ge h\right\}\]
be the first time that $h$ is compensated by paying dividends with strategy $\pi$. The strategy $\yw{\pi}$ does not pay dividends before time $\t^h$.
 We define the strategy $\yw{\pi}$ as follows:
\[L_t^{\yw{\pi}}=\zb(L_t^\pi-h\yb)^+\mathbf{1}_{\zb\{\t^h\le t\le T\yb\}}\mathbf{1}_{\zb\{t<\yw{\t}\yb\}}+X_T^{\yw{\pi},x-h}\mathbf{1}_{\zb\{t=T<\t^h\qx\yw{\t}\yb\}},\]
where $(L_t^\pi-h)^+=\max\zb\{0,L_t^\pi-h\yb\}$ and $\yw{\tau}=\inf\{t\ge s: X_t^{\yw{\pi},x-h}<0\}$
is the ruin time of $X_t^{\yw{\pi},x-h}$.
Note that $\yw{\pi}$ is admissible. 
From the definition of $\yw{\pi}$, we obtain
\begin{align}\label{a33}
&\mathbb{E}_{sxw}\left[\jf_{s}^{(\t^\pi\qx T)+}e^{-c(t-s)}dL_t^\pi-\jf_{s}^{(\yw{\t}\qx T)+}e^{-c(t-s)}dL_t^{\yw{\pi}}:S\le \yw{\t}\right]\le L_{(S\qx T)+}^\pi\!-\!L_{(S\qx T)+}^{\yw{\pi}}\le h.\\\label{a34}
&\mathbb{E}_{sxw}\left[\jf_{s}^{(\t^\pi\qx T)+}e^{-c(t-s)}dL_t^\pi-\jf_{s}^{(\yw{\t}\qx T)+}e^{-c(t-s)}dL_t^{\yw{\pi}}:S>\yw{\t}, \yw{\t}=\t^\pi\right]\le h.
\end{align}
As for any strategy $\pi$, the corresponding cost function is bounded, we obtain
\begin{align}\bbh
&\mathbb{E}_{sxw}\left[\jf_{s}^{(\t^\pi\qx T)+}e^{-c(t-s)}dL_t^\pi-\jf_{s}^{(\yw{\t}\qx T)+}e^{-c(t-s)}dL_t^{\yw{\pi}}:\t^h>\yw{\t},\yw{\t}\bdy\t^\pi\right]\\\bbh\le&
\left(x+pT\right)\mathbb{P}(\t^h>\yw{\t},\yw{\t}\bdy\t^\pi)
\le \left(x+pT\right)\mathbb{P}\left\{\Delta X_{\yw{\t}}^{\yw{\pi},x-h}\in (X_{\yw{\t}-}^{\yw{\pi},x-h}, X_{\yw{\t}-}^{\yw{\pi},x-h}+h)\right\}\\\bbh
\le &\left(x+pT\right)\jf_0^\wq \mathbb{P}\left\{\Delta X_{\yw{\t}}^{\yw{\pi},x-h}\in (y, y+h)|X_{\yw{\t}-}^{\yw{\pi},x-h}=y\right\}F_{X_{\yw{\t}-}^{\yw{\pi},x-h}}(dy)\\\label{a32}
=&\left(x+pT\right)\jf_0^\wq [G(y+h)-G(y)]F_{X_{\yw{\t}-}^{\yw{\pi},x-h}}(dy),
\end{align}
where $G$ denotes the common distribution function of the claim size $U_i$'s. As $G$ is uniformly continuous, we observe that for any $\e>0$, there exists a constant $\d_1(x)>0$ such that for all $h\le \d_1(x)$, we have $G(y+h)-G(y)\le \fs{\e}{2(x+pT)}$. Combining this with (\ref{a33}), (\ref{a34}) and (\ref{a32}), we obtain that for  any $\e>0$, there exists a $\d(x)=\min\{\fs{\e}{2}, \d_1(x)\}>0$ such that for all $h<\d(x)$ and  any strategy $\pi\in U_{ad}^{x,w}[s,T]$, it holds that
\[
J(s,x,w;\pi)\le J(s,x-h,w;\yw{\pi})+\e\le V(s,x-h,w)+\e.
\]
As $\pi\in U_{ad}^{x,w}[s,T]$ is arbitrary, we obtain that for a sufficiently small $h$, we have $V(s,x,w)-V(s,x-h,w)\le \e$. This shows that $V$ is continuous with respect to $x$. Thus, for all compact sets, $V$ is uniformly continuous with respect to $x$.  $\hfill\Box$
\begin{prop}\label{a31} Under assumption \ref{a52}, for any $h>0$ and $0\le s<s+h\le T$, it holds that
(1)$V(s+h,x,w+h)-V(s,x,w)\le \left[1-\exp\zb\{-ch-\int_w^{w+h}\lambda (u)du\yb\}\right]V(s+h,x,w+h)$.\\
(2)$V(s,x,w+h)-V(s,x,w)\le 2ph+\left[1-\exp\zb\{-ch-\int_w^{w+h}\lambda (u)du\yb\}\right]V(s+h,x,w+h)$.
\end{prop}
$\mathbf{Proof}$\q
(1) For all strategies $\pi\in U_{ad}^{s,w+h}[s+h,T]$, $\yw{\pi}^h\in U_{ad}^{x,w}[s,T]$  is defined so that on set $\{T_1^{s,w}>h\}$ the insurance  company pays  dividends at  rate $p$ during time $[s,s+h)$, and then follows the strategy $\pi$ during $t\in [s+h,T]$; on set $\{T_1^{s,w}\le h\}$, $\yw{\pi}^h$ pays no  dividends before time $T$ and then distributes the entire wealth as dividends at time $T$ if $X_T^{\yw{\pi}^h}>0$. The rest of the proof is exactly the same as the proof of proposition 5.1 of \cite{bai2017optimal}. Depending on whether there are claims during the time interval $[s,s+h]$, we can deduce that
\bea\label{chanteur}
V(s,x,w)\ge \exp\left\{-ch-\int_w^{w+h}\lambda(u)du\right\} V(s+h,x,w+h).
\eea
(2) For the second part of the proof, we can calculate the following directly:
 \begin{eqnarray}\label{a29}
V(s,x,w\!+\!h)\!-\!V(s,x,w)\!=\!V(s,x,w\!+\!h)\!-\!\!V(s\!+\!h,x,w\!+\!h)\!+\!V(s\!+\!h,x,w\!+\!h)\!-\!V(s,x,w).
\end{eqnarray}
Combining (\ref{a28}) and (\ref{chanteur}) with (\ref{a29}), we obtain
\[V(s,x,w+h)-V(s,x,w)\le 2ph\!+\!\left\{1\!-\!\exp\left[-ch\!-\!\int_w^{w+h}\lambda(u)du\right]\right\}V(s\!+\!h,x,w\!+\!h).\]
This completes the proof. $\hfill\Box$
\begin{prop} Under assumption \ref{a52}, the value function has the following property.
\[\lim_{h\downarrow 0}[V(s,x,w)-V(s,x,w+h)]= 0,\]
uniformly in a compact set ${\mathfrak{D}}$.
\end{prop}
$\mathbf{Proof}$\q From  Proposition \ref{a31}$-(2)$, we note  that $\lim_{h\downarrow 0}[V(s,x,w)-V(s,x,w+h)]\ge 0$.
We only need to prove the opposite inequality. By a direct calculation and Proposition \ref{a50}-(1), we obtain
\begin{eqnarray}
\nonumber
V(s,x,w)\!-\!V(s,x,w+h)\!\!\!&=&\!\!\!V(s,x,w)\!\!-\!V(s\!+\!h,x,w\!+\!h)\!+\!V(s\!+\!h,x,w\!+\!h)\!-\!V(s,x,w\!+\!h)\\\label{6}\!&\le&\!V(s,x,w)\!\!-\!\!V(s+h,x,w+h).
\end{eqnarray}
For any strategy $\pi\in U_{ad}^{x,w}[s,T]$, depending on whether there are claims between $[s,s+h]$, we obtain
\begin{align}
\nonumber
&\mathbb{E}_{sxw}\left[\int_{s}^{(\tau^{\pi}\land T)+}e^{-c(t-s)}dL_t^\pi\right]-V(s+h,x,w+h)\\\nonumber=&\mathbb{E}_{sxw}\left[\int_{s}^{(\tau^{\pi}\land T)+}e^{-c(t-s)}dL_t^\pi|T_1^{s,w}\le h\right]\mathbb{P}(T_1^{s,w}\le h)\\\label{a36}&+\mathbb{E}_{sxw}\left[\int_{s}^{(\tau^{\pi}\land T)+}e^{-c(t-s)}dL_t^\pi|T_1^{s,w}>h\right]\mathbb{P}(T_1^{s,w}>h)-V(s+h,x,w+h).
\end{align}
For the first term of the right side of (\ref{a36}), because
$
\mathbb{P}(T_1^{s,w}\le h)=1-\exp[-\jf_w^{w+h}\l (u)du],
$
we know that for a sufficiently small $h$, there exists a constant $Q_1(\mathfrak{D})>0$ such that
\bea\label{a44}
\mathbb{E}_{sxw}\left[\int_{s}^{(\tau^{\pi}\land T)+}e^{-c(t-s)}dL_t^\pi|T_1^{s,w}\le h\right]\mathbb{P}(T_1^{s,w}\le h)\le Q_1(\mathfrak{D})h,
\eea
which means that for all $\e>0$, there exists a constant $\d_1>0$, such that for all $0<h<\d_1$ and all strategies $\pi\in U_{ad}^{x,w}[s,T]$, it holds that
\bea\label{8}
\mathbb{E}_{sxw}\left[\int_{s}^{(\tau^{\pi}\land T)+}e^{-c(t-s)}dL_t^\pi|T_1^{s,w}\le h\right]\mathbb{P}(T_1^{s,w}\le h)\le\fs{\e}{3}.
\eea
We estimate the second and third terms of the right side of (\ref{a36}) as follows:
\begin{align}\bbh
&\mathbb{E}_{sxw}\left[\int_{s}^{(\tau^{\pi}\land T)+}e^{-c(t-s)}dL_t^\pi|T_1^{s,w}>h\right]\mathbb{P}(T_1^{s,w}>h)-V(s+h,x,w+h)\\\label{320}
\le&\mathbb{E}_{sxw}\left[\int_{s}^{(\tau^{\pi}\land T)+}e^{-c(t-s)}dL_t^\pi|T_1^{s,w}>h\right]-V(s+h,x,w+h).
\end{align}
As on $\{T_1^{s,w}>h\}$, there is no claim on $[s,s+h]$, and it is clear that $\t^\pi>s+h$. Thus, we can deduce that
\begin{align}\bbh
\mathbb{E}_{sxw}\!\left[\!\int_{s}^{(\tau^{\pi}\land T)+}\!\!\!e^{-c(t-s)}dL_t^\pi|T_1^{s,w}\!\!>\!\!h\right]
&\!=\!\!\mathbb{E}_{sxw}\!\left[\int_{s}^{s+h}\!\!e^{-c(t-s)}dL_t^\pi\!\!+\!\!\int_{(s+h)-}^{(\tau^{\pi}\land T)+}\!\!e^{-c(t-s)}dL_t|T_1^{s,w}\!\!>\!\!h\right]
\\\bbh
&\!\le\! \mathbb{E}_{sxw}\left[L_{s+h}^\pi\!\!+\!\!V(s\!+\!h,X_{s+h}^{\pi},w\!\!+\!\!h)|T_1^{s,w}\!\!>\!\!h\right]\\\label{a35}
&\!=\!\mathbb{E}_{sxw}\!\left[x\!\!+\!\!ph\!-\!X_{s+h}^{\pi}\!+\!V(s\!\!+\!\!h,X_{s+h}^{\pi},w\!\!+\!\!h)|T_1^{s,w}\!\!>\!\!h\right].
\end{align}
The last equality sign of (\ref{a35}) is based on the fact that $L_{s+h}^\pi=x+ph-X_{s+h}^\pi$ on $\{T_1^{s,w}\ge h\}$. Then, we can obtain 
\begin{align}\bbh
&\mathbb{E}_{sxw}\left[\int_{s}^{(\tau^{\pi}\land T)+}e^{-c(t-s)}dL_t|T_1^{s,w}>h\right]-V(s+h,x,w+h)\\\label{a37}
\le &\mathbb{E}_{sxw}\left[x+ph-X_{s+h}^{\pi}+V(s+h,X_{s+h}^{\pi},w+h)-V(s+h,x,w+h)|T_1^{s,w}>h\right].
\end{align}
As $0\le X_{s+h}^\pi\le x+ph$ on $\{T_1^{s,w}>h\}$, conditioning  on whether  $X_{s+h}^\pi$ belongs to $[0,x]$, we can estimate (\ref{a37}) as follows:
\begin{eqnarray}\nonumber & &
\mathbb{E}_{sxw}\left[x+ph-X_{s+h}^{\pi}+V(s+h,X_{s+h}^{\pi},w+h)-V(s+h,x,w+h)|T_1^{s,w}>h\right]\\\nonumber &=&\!\!
\mathbb{E}_{sxw}[x\!+\!ph\!-\!X_{s+h}^{\pi}\!\!+\!\!V(s\!+\!h,X_{s+h}^{\pi},w\!+\!h)\!\!-\!\!V(s\!+\!h,x,w\!+\!h)|0\!\le\! X_{s+h}^\pi\!\le\!x,T_1^{s,w}>h]\\\bbh&&\times\mathbb{P}(0\le X_{s+h}^\pi\le x|T_1^{s,w}>h)\!+\!\mathbb{P}(x\!<\!X_{s+h}^\pi\le x+ph|T_1^{s,w}\!>\!h)\!\times\!\mathbb{E}_{sxw}[x\!+\!ph\!-\!X_{s+h}^{\pi}\!\\\label{a41} & & +\!V(s\!+\!h,X_{s+h}^{\pi},w\!+\!h)\!\!-\!\!V(s\!+\!h,x,w\!+\!h)|x\!<\!X_{s+h}^\pi\!\le \!x\!+\!ph, T_1^{s,w}\!>\!h].
\end{eqnarray}
As on  $\{0\le X_{s+h}^\pi\le x,T_1^{s,w}>h\}$,
\bea\label{a40}
V(s+h,X_{s+h}^{\pi},w+h)-V(s+h,x,w+h)\le-(x-X_{s+h}^\pi).  \eea
On $\{x<X_{s+h}^\pi\le x+ph,T_1^{s,w}>h\}$,
\begin{align}\label{a39}
V(s\!+\!h,X_{s+h}^\pi,w\!+\!h)\!-\!V(s\!+\!h,x,w+h)\!\le\! V(s\!+\!h,x\!+\!ph,w\!+\!h)\!-\!V(s\!+\!h,x,w\!+\!h).
\end{align}
Substituting (\ref{a40}) and (\ref{a39}) into (\ref{a41}), we deduce that
\begin{align}\bbh
&\mathbb{E}_{sxw}\left[x+ph-X_{s+h}^{\pi}+V(s+h,X_{s+h}^{\pi},w+h)-V(s+h,x,w+h)|T_1^{s,w}>h\right]\\\label{a42}
\le\;& ph+V(s+h,x+ph,w+h)-V(s+h,x,w+h).
\end{align}
Substituting (\ref{a42}) into (\ref{a37}) and combining this with (\ref{320}), we obtain
\begin{align}\bbh
&\mathbb{E}_{sxw}\left[\int_{s}^{(\tau^{\pi}\land T)+}e^{-c(t-s)}dL_t^\pi|T_1^{s,w}>h\right]\mathbb{P}(T_1^{s,w}>h)-V(s+h,x,w+h)\\\label{a43}
\le\;& ph\!+\!V(s\!+\!h,x\!+\!ph,w\!+\!h)\!\!-\!\!V(s\!+\!h,x,w\!+\!h).
\end{align}
By Proposition \ref{a38}, we obtain
\[ \lim_{h\downarrow 0}V(s+h,x+ph,w+h)-V(s+h,x,w+h)=0.\]
Thus, for all $\e>0$, there exists a constant $\d_2>0$ such that for all $0<h<\d_2$, it holds that
\[V(s+h,x+ph,w+h)-V(s+h,x,w+h)\le\fs{\e}{3}.\]
Setting $\d_3=\min\{\d_2,\fs{\e}{3p}\}$,  for all $0<h<\d_3$, from (\ref{a43}), we obtain
\bea\label{7}
\mathbb{E}_{sxw}\left[\int_{s}^{(\tau^{\pi}\land T)+}e^{-c(t-s)}dL_t^\pi|T_1^{s,w}>h\right]\mathbb{P}(T_1^{s,w}>h)-V(s+h,x,w+h)\le\fs{2\e}{3}.
\eea
Combining (\ref{a36}) and (\ref{8}) with (\ref{7}), we know that
for all $\e>0$, there exists a constant  $\d_4=\min\{\d_1,\d_3\}$ such that  for all $0<h<\d_4$ and all strategies $\pi\in U_{ad}^{x,w}[s,T]$, it holds that
\[\mathbb{E}_{sxw}\left[\int_{s}^{(\tau^{\pi}\land T)+}e^{-c(t-s)}dL_t^\pi\right]-V(s+h,x,w+h)\le\e.\]
As $\pi\in U_{ad}^{x,w}[s,T]$ is arbitrary, we have
\bea\label{9}
V(s,x,w)-V(s+h,x,w+h)\le \e.
\eea
From (\ref{6}) and (\ref{9}), we obtain
\[\lim_{h\downarrow 0}[V(s,x,w)-V(s,x,w+h)]\le 0.\]
This completes the proof. $\hfill\Box$
\begin{cor}
Combining all the above properties, we observe that $V$ is continuous on $D$ and uniformly continuous in any compact set $\mathfrak{D}$.
\end{cor}

\section{Dynamic Programming Principle}
In this section, we develop the dynamic programming principle (DPP) for our optimization problem. We start by proving an important lemma.
{{\begin{lem}\label{a4}
For any compact set $\mathfrak{D}$ and $\e >0$,
 there exists a constant $\delta>0$ independent of $(s,x,w)\in\mathfrak{D}$, such that for any $\pi\in U_{ad}^{x,w}[s,T]$ and $0<h<\delta$, we can find $\check{\pi}\in{U}_{ad}^{x,w-h}[s,T]$ such that
\[
J(s,x,w;\pi)-J(s,x,w-h;\check{\pi})\leq\e, \quad\quad \forall\; (s,x,w)\in \mathfrak{D}.
\]
\end{lem}
$\mathbf{Proof}$\q
Let $\pi\in U_{ad}^{x,w}[s,T]$. For any $h>0$, we consider the following two strategies.
\begin{enumerate}
\item Define $\pi_1\in U_{ad}^{x,w-h}[s-h,T]$ such that in the case of  $\{T_1^{s-h,w-h}>h\}$, $\pi_1$ pays dividends at rate $p$ during the time interval $[s-h,s)$ and then follows strategy $\pi$ during the time interval $t\in[s,T]$; in the case of $\{T_1^{s-h,w-h}<h\},\pi_1$ pays a jump dividend $x$ at time $s-h$ and then pays dividends at rate $p$ until the ruin time.
\item Define $\check{{\pi}}\in U_{ad}^{x,w-h}[s,T]$ such that \[\check{{\pi}}(t)={\pi_1}(t-h)\mathbf{1}_{\{t<T\}}+\pi^{all}\mathbf{1}_{t=T}\mathbf{1}_{X_{T}^{\check{\pi}}>0},\]
 where $\pi^{all}$ denotes the strategy of paying the entire wealth as a jump  dividend.
\end{enumerate}
We can calculate directly that
\[
J(s,x,w;\pi)\!\!-\!\!J(s,x,w\!\!-\!\!h;\check{\pi})\!\!=\!\!J(s,x,w;\pi)\!\!-\!\!J(s\!\!-\!\!h,x,w\!\!-\!\!h;\pi_1)\!\!+\!\!J(s\!\!-\!\!h,x,w\!\!-\!\!h;\pi_1)\!\!-\!\!J(s,x,w\!\!-\!\!h;\check{\pi}).
\]
By the proof of proposition \ref{a31}, we obtain
\begin{align}\bbh
 J_1\!&\dy\! J(s,x,w;\!\pi)\!\!-\!\!J(s\!-\!h,x,w\!-\!h;\pi_1)
\!\! \le\!\! J(s,x,w;\pi)\zb[1\!-\!\exp\zb\{\!-\!ch\!-\!\jf_{w\!-\!h}^{w}\!\!\l(u)du\yb\}\yb]\\\bbh
 &\le(x+pT)\zb[1-\exp\zb\{-ch-\jf_{w\!-\!h}^{w}\l(u)du\yb\}\yb].
\end{align}
Using the fact that $\lim_{h\yjt0}\zb[1-\exp\zb\{-ch-\jf_{w\!-\!h}^{w}\l(u)du\yb\}\yb]=0$, we observe that for all $\e>0$, there exists a constant $\d_1>0$ such that for all $0<h<\d_1$, it holds that $J_1<\fs{\e}{2}$. From the proof of proposition \ref{a50}-(2)(e.g., the continuity of the value function on $s$), we obtain
\[
J_2\dy J(s-h,x,w-h;\pi_1)-J(s,x,w-h;\check{\pi})\le ph.
\]
Thus, for all $\e>0$, there exists a constant $\d_2>0$ such that for all $h<\d_2$, there exists $\check{\pi}\in U_{ad}^{x,w-h}[s,T]$ such that $J_2<\fs{\e}{2}$. Taking $\d=\min\{\d_1,\d_2\}$, we obtain the desired conclusion.
 $\hfill\Box$}}

\no
For any strategy $\pi\in{U}_{ad}^{x,w}(s,T)$,  denote $R_t^\pi$ $=$ $R_t^{\pi,s,x,w}$\;$=(t,X_t^{\pi,s,x,w},W_t^{s,w})$ for simplicity.
\begin{thm} Assume that assumption~\ref{a52} holds; then,
for any $(s,x,w)\in D$ and any stopping time $\tau\in[s,T]$, it holds that
\[
V(s,x,w)=\sup_{\pi\in{U}_{ad}^{x,w}[s,T]}\mathbb{E}_{sxw}\left[\int_{s}^{\tau\wedge \tau^{\pi}}e^{-c(t-s)}dL_t^\pi\!+\!e^{-c(\tau\wedge\tau^{\pi}-s)}V(R_{\tau\wedge\tau^\pi}^\pi)\right].
\]
\end{thm}
$\mathbf{Proof}$ \q First, we  prove  that the above equation holds for any deterministic $\tau=s+h$, where $h\in(0,T-s)$. In other words, we prove that
\[
V(s,x,w)=v(s,x,w;s+h),
\]
where
\[v(s,x,w;s+h)\!\triangleq\!\sup_{\pi\in{U}_{ad}^{x,w}[s,T]}\mathbb{E}_{sxw}\left[\int_{s}^{(s+h)\wedge\tau^\pi}\!\!\!e^{-c(t-s)}dL_t^\pi\!+\!e^{-c\zb((s+h)\wedge \tau^\pi-s\yb)}V(R_{(s+h)\wedge \tau^\pi}^\pi)\right].\]
First, we show that $V(s,x,w)\le v(s,x,w;s+h)$. The proof of this statement is exactly the same as the proof of Theorem 6.2 in \cite{bai2017optimal}. We omit this proof here.

In what follows, we prove that $V(s,x,w)\ge v(s,x,w;s+h)$. Denote $x_{max}\dy x+pT$. For a compact set $[0,x_{\max}]\times[0,T]^2$, due to the uniform continuity of $V$ on $(x,w)$, we know that there exists a constant $\d_3>0$ such that for all $0<x_A-x_B<\d_3$ and $0<w_A-w_B<\d_3$, it holds that  $|V(s,x_A,w_A)-V(s,x_B,w_B)|<\fs{\e}{2}$.
 Let $\d_4>0$ be the constant in Lemma \ref{a4}. Denote $\d\dy\min\{\d_3,\d_4\}$.  Let $0=x_0<x_1<x_2<\cdots<x_n=x_{\max}$ and $0=w_0<w_1<w_2<\cdots<w_n=T$ be a partition of $[0,x_{\max}]\times [0,T]$ such that $x_{i+1}-x_i<\d$ and $w_{j+1}-w_j<\delta$. Take $D_{ij}=[x_{i-1},x_i)\times[w_{j-1},w_j)$, $i,j\in\mathbb{N}.$ For $0\le s<s+h<T$ and $0\le i,j\le n$, we choose $\pi^{i,j}\in{U}_{ad}^{x_i,w_j}[s+h,T]$ such that
\begin{eqnarray}\label{10}
J(s+h,x_{i},w_j;\pi^{i,j})>V(s+h,x_{i},w_j)-\e.
\end{eqnarray}
For each $x\in [x_{i-1},x_i)$, we define $\yw{\pi}^{i,j}\in U_{ad}^{x,w_j}[s+h,T]$ such that $\yw{\pi}^{i,j}$ pays a jump dividend $x-x_{i-1}$ and then follows the strategy $\pi^{i-1,j}$, meaning that
\bea\label{2}
J(s+h,x,w_j;\yw{\pi}^{i,j})=J\zb(s+h,x_{i-1},w_j;\pi^{i-1,j}\yb)+x-x_{i-1}.
\eea
For each $w\in [w_{j-1},w_j),x\in[x_{i-1},x_i)$ and $\yw{\pi}^{i,j}\in U_{ad}^{x,w_j}[s+h,T]$, by Lemma \ref{a4}, we can define strategy $\hat{\pi}^{i,j}\in U_{ad}^{x,w}[s+h,T]$ such that
\bea\label{3}
J(s+h,x,w;\hat{\pi}^{i,j})\ge J(s+h,x,w_j;\yw{\pi}^{i,j})-\e.
\eea
Thus, combining (\ref{10}), (\ref{2}) and (\ref{3}) we obtain that for each $(x,w)\in D_{ij}$,
\begin{align}\bbh
J(s+h,x,w;\yj{\pi}^{i,j})&\ge J(s+h,x,w_j;\yw{\pi}^{i,j})-\e=J(s+h,x_{i-1},w_j;\pi^{i-1,j})+x-x_{i-1}-\e\\\bbh
&\ge V(s+h,x_{i-1},w_j)-2\e\ge V(s+h,x,w)-3\e.
\end{align}
The last inequality holds because $V$ is uniformly continuous with respect to $x$ and $w$ on a compact set.
For any strategy $\pi\in{U}_{ad}^{x,w}[s,T]$, define the new strategy $\pi^*$ as follows.
\[
\pi^*(t)=\pi(t) \textbf{1}_{[s,s+h)}(t)+\sum_{i=0}^{n-1}\sum_{j=0}^{n-1}\textbf{1}_{D_{ij}}(X_{s+h}^\pi,W_{s+h})\textbf{1}_{[s+h,T]}(t)\hat{\pi}^{i,j}(t).
\]
We can verify that $\pi^*\in{U}_{ad}^{x,w}[s,T]$  and $\{\tau^{\pi^*}\le s+h\}=\{\tau^\pi\le s+h\}$. If $\tau^{\pi^*}> s+h$, we have
\[
J(s+h,X_{s+h}^\pi,W_{s+h};\pi^*)\ge V(s+h, X_{s+h}^\pi,W_{s+h})-3\e, \q\mathbb{P}-a.s. \;\bh{on}\; \{\tau^{\pi^*}> s+h\}.
\]
Consequently,
\begin{align}\nonumber
V(s,x,w)\!&\ge \!\! J(s,x,w;\pi^*)\!=\!\mathbb{E}_{sxw}\!\!\left[\int_{s}^{(s+h)\wedge \tau^\pi}\!\!\!\!e^{-c(t-s)}dL_t\!+\!\textbf{1}_{{\{\tau^\pi>s+h}\}}\!e^{-ch}\!\!\int_{s}^{\tau^{\pi^*}}\!\!\!\!\!e^{-c(t-(s+h))}dL_t^*\right]\\\label{12}
&\ge\mathbb{E}_{sxw}\zb[\int_{s}^{(s+h)\wedge \tau^\pi}e^{-c(t-s)}dL_t+e^{-c((s+h)\qx\t^\pi-s)}V(R_{(s+h)\qx\t^\pi}^\pi)\yb]-3\e.
\end{align}
In the last inequality, we used the fact that $\textbf{1}_{{\{\tau^\pi\le s+h}\}}V(R_{(s+h)\qx\t^\pi}^\pi)=0.$ As $\pi\in{U}_{ad}^{x,w}[s,T]$ is arbitrary, from (\ref{12}) we observe that $V(s,x,w)\ge v(s,x,w;s+h)-3\e$. As $\e>0$ is arbitrary, we obtain that $V(s,x,w)\ge v(s,x,w;s+h)$ for all deterministic  $\tau=s+h$. Since the proof of $\t$ being a stopping time is exactly the same as the proof of Theorem 6.2 in  \cite{bai2017optimal}, we omit the rest of the proof.  $\hfill\Box$
\section{The Hamilton-Jacobi-Bellman equation}
The corresponding HJB equation of the  problem is
\begin{eqnarray}
\begin{cases}
\label{1}
\max\left\{1-V_x,\hx{L}[V]\right\}(s,x,w)=0, (s,x,w)\in\hx{D};\\
V(T,x,w)=x,
\end{cases}
\end{eqnarray}
where $\hx{L}[\cdot]$ is a first-order integro-differential operator for $\varphi\in \mathbb{C}^{1,1,1}(D)$, where $\mathbb{C}^{1,1,1}(D)$ denotes the set of all continuously differentiable  functions on $D$, and
\begin{eqnarray*}
\hx{L}[\varphi](s,x,w)\dy[-c\varphi\!+\!\varphi_s\!+\!p\varphi_x+\varphi_w](s,x,w)+\lambda(w)\left[\int_0^x\varphi(s,x-u,0)dG(u)-\varphi(s,x,w)\right].
\end{eqnarray*}
From now on, we denote the integral part of (\ref{1}) as $I[\vf]$ for simplicity, i.e.,
\[
I[\vf](s,x,w)=\int_0^x\varphi(s,x-u,0)dG(u)-\varphi(s,x,w).
\]
While we want to explore more regularity properties of the value function, unfortunately, in many applications, the value function $V(s,x,w)$ is not necessarily smooth, or it can be very difficult to prove its differentiability. Therefore, we need to introduce the notion of weak solutions, namely, viscosity solutions.

We recall that the notion of viscosity solutions was introduced by Crandall and Lions \cite{crandall1983viscosity} for first-order equations and Lions \cite{Lions1983,Lions1983Optimalone} for second-order equations. The notion of viscosity solutions of integro-differential equations was pursued by Soner\cite{soner1986optimal}. The viscosity solution concept for  fully nonlinear partial differential equations has been proven to be extremely useful in control theory because it does not need the differentiability of the value function.
 Instead, it merely requires continuity of the value function to define the viscosity solution. We refer to the user's guide of Crandall, Ishii and Lions\cite{crandall1992user} for an overview of the theory of viscosity solutions and their applications. Using the notion of a viscosity solution, we   prove that the value function is the (viscosity) solution of the corresponding equation.  The viscosity solution approach is becoming a well-established approach to studying stochastic control problem; see, e.g., the books \cite{fleming2006controlled,yong1999stochastic}.
Now, we provide the definition of the constrained viscosity solution of the HJB equation (\ref{1}).
\begin{defn}
Let $O\subseteq\hx{D}^*$ be a subset such that  $\partial_T O\dy\{(T,y,v)\in \partial O\}\neq \emptyset$, where $\bar{O}$ is the closure of $O$. Denote USC$(O)$ as the set of all upper semicontinuous functions on $O$ and $LSC(O)$ as the set of all lower semicontinuous functions on $O$.

(a) Let $v\in \bh{USC}(O)$; we call $v$ a viscosity subsolution of (\ref{1}) on $O$ if $v(T,y,v)\le y$ for $(T,y,v)\in \partial_T O$; for any $(s,x,w)\in O$ and $\vf \in\mathbb{C}^{1,1,1}(\bar{O})$ such that $[v-\vf](s,x,w)=\max_{(t,y,v)\in O}[v-\vf](t,y,v)$, it holds that
\[
\max\{1-\vf_x,\hx{L}[\vf]\}(s,x,w)\ge 0.
\]

(b) Let $v\in \bh{LSC}(O)$; we call $v$ a viscosity supersolution of (\ref{1}) on O if $v(T,y,v)\ge y$ for all $(T,y,v)\in \partial_T O $; for any $(s,x,w)\in O$ and $\vf \in\mathbb{C}^{1,1,1}(\bar{O})$ such that $0=[v-\vf](s,x,w)=\min_{(t,y,v)\in O}[v-\vf](t,y,v)$, it holds that
\[
\max\{1-\vf_x,\hx{L}[\vf]\}(s,x,w)\le 0.
\]
In particular, we call $u$ a ``constrained viscosity solution'' of (\ref{1}) on $\hx{D}^*$ if it is both a viscosity subsolution on $\hx{D}^*$ and a viscosity supersolution on $\hx{D}.$
\end{defn}
%
There is an equivalent formulation of viscosity solutions; the proof of equivalence of definitions is standard (e.g., see\cite{awatif1991equqtions,benth2001optimal}). In this paper, we use both definitions interchangeably.

Given a continuously differentiable function $\vf$ and a continuous function $u$, we define the operator
\begin{align}\bbh
\hx{L}[u,\vf](s,x,w)\!=\![-cu\!+\!\varphi_s\!+\!p\varphi_x\!+\!\varphi_w](s,x,w)\!+\!\lambda(w)\left[\int_0^xu(s,x\!-\!\a,0)dG(\a)-u(s,x,w)\right].
\end{align}
\begin{defn}
Let $v\in \bh{USC}(O)$; we call $v$ a viscosity subsolution of (\ref{1}) on $O$ if $v(T,y,v)\le y$ for $(T,y,v)\in \partial_T O$; for any $(s,x,w)\in O$, $\vf \in\mathbb{C}^{1,1,1}(\bar{O})$ such that $0=[v-\vf](s,x,w)=\max_{(t,y,v)\in O}[v-\vf](t,y,v)$, it holds that
\[
\max\{1-\vf_x,\hx{L}[v,\vf]\}(s,x,w)\ge 0.
\]

Let $v\in \bh{LSC}(O)$; we call $v$ a viscosity supersolution of (\ref{1}) on O if $v(T,y,v)\ge y$ for all $(T,y,v)\in \partial_T O $; for any $(s,x,w)\in O$ and $\vf \in\mathbb{C}^{1,1,1}(\bar{O})$ such that $0=[v-\vf](s,x,w)=\min_{(t,y,v)\in O}[v-\vf](t,y,v)$, it holds that
\[
\max\{1-\vf_x,\hx{L}[v,\vf]\}(s,x,w)\le 0.
\] 
\end{defn}
\begin{thm}
Assume that assumption \ref{a52} holds. Then,
the value function is a constrained viscosity solution of the HJB equation (\ref{1}) on $\hx{D}^*$.
\end{thm}
$\mathbf{Proof}$\q
$Supersolution$ \;Given $(s,x,w)\in\hx{D}$, let $\varphi\in{\mathbb{C}}^{1,1,1}(\hx{D})$ be such that $V-\varphi$ reaches its minimum at $(s,x,w)$ with $V(s,x,w)=\varphi(s,x,w)$. Consider the strategy $\pi^0$ with dividend rate $l_0$; $T_1^{s,w}$ denotes the time of the first claim. Denote $\tau_s^h\dy s+h\wedge T_1^{s,w}, R_t\dy(t,X_t^{\pi^0,s,x,w},W_t^{s,w})$. By the dynamic programming principle, we have
\begin{align}\nonumber
V(s,x,w)&=\sup_{\pi\in{U}_{ad}^{x,w}[s,T]}\mathbb{E}_{sxw}\left[\int_{s}^{\tau_s^h}e^{-c(t-s)}dL_t^\pi+e^{-c(\tau_s^h-s)}
V(R_{\tau_s^h}^\pi)\right]
\\\nonumber
&\ge \mathbb{E}_{sxw}\left[\int_{s}^{\tau_s^h}e^{-c(t-s)}l_0dt+e^{-c(\tau_s^h-s)}V(R_{\tau_s^h}^\pi)\right]\\\nonumber
&=\mathbb{E}_{sxw}\left[\int_{s}^{\tau_s^h}\!\!e^{-c(t-s)}l_0dt\right]\!+\!\mathbb{E}_{sxw}\left\{e^{-c(\tau_s^h-s)}\left[V(R_{\tau_s^h})-V (R_{{\tau_s^h}-})\right]\mathbf{1}_{\{T_1^{s,w}<h\}}\right\}\\\nonumber &\q +\mathbb{E}_{sxw}\left[e^{-c(\tau_s^h-s)}\varphi (R_{{\tau_s^h}-})\right].
\end{align}
Using the fact that $V(s,x,w)=\varphi(s,x,w)$, we obtain
\begin{align}\nonumber
0&\ge\mathbb{E}_{sxw}\left[\int_{s}^{\tau_s^h}e^{-c(t-s)}l_0dt\right]+\mathbb{E}_{sxw}\left\{e^{-c(\tau_s^h-s)}\zb[V(R_{\tau_s^h})-V (R_{{\tau_s^h}-})\yb]\textsl{1}_{\{T_1^{s,w}<h\}}\right\}\\\label{5} &\quad  +\mathbb{E}_{sxw}\left[e^{-c(\tau_s^h-s)}\varphi (R_{{\tau_s^h}-})-\varphi(s,x,w)\right]\triangleq I_1+I_2+I_3.
\end{align}
where $I_i,i=1,2,3$ are the three terms on the right side above. Clearly, we have
\begin{align}\bbh
I_1&=l_0\mathbb{E}_{sxw}\left[\int_{s}^{s+h\wedge T_1^{s,w}}e^{-c(t-s)}dt:T_1^{s,w}\le h\right]
+l_0\mathbb{E}_{sxw}\left[\int_{s}^{s+h}e^{-c(t-s)}dt:T_1^{s,w}>h\right]\\\label{I1} &= l_0\int_0^h\!\!\int_{s}^{s+z}\!e^{-c(t-s)}\lambda(w\!+\!z)e^{-\int_w^{w+z}\lambda(u)du}dtdz\!+\!\!l_0\!\left[1\!-\!e^{-\int_w^{w+h}\lambda (u)du}\right]\!\!\int_{s}^{s+h}e^{-c(t-s)}dt.
\end{align}
As $\t_s^h=s+T_1^{s,w}$ on $\{T_1^{s,w}<h\}$, we have
\begin{align}\label{I2}
I_2\!=\!\mathbb{E}_{sxw}\!\zb[\int_0^\infty\!\! \int_0^h\!e^{-ct}\left[V (s+t,X_{s+t-}^{\pi^0}\!-u,0)\!-\!V(s+t,X_{{s+t}-}^{\pi^0},W_{{s+t}-}^{s,w})\right]dF_{T_1^{s,w}}(t)dG(u)\yb].
\end{align}
As there are no jumps on $[s,\tau_s^h)$, using It$\hat{\bh{o}}$'s formula, we obtain
\begin{align}\nonumber
I_3&=\mathbb{E}_{sxw}\left[\int_{s}^{\tau_s^h}e^{-c(u-s)}\left[-c\varphi +\varphi_t-l_0\varphi_x+\varphi_xp+\varphi_w\right](R_u^0)du\right]\\\label{I3}&= \mathbb{E}_{sxw}\zb[\int_{s}^{s+h}\bar{F}_{T_1^{s,w}}(u-s)e^{-c(u-s)}\left[-c\varphi\!+\!\varphi_t-\varphi_xl_0+\varphi_xp+\varphi_w\right](R_u^0)du\yb].
\end{align}
Recall that $F_{T_1^{s,w}}(t)=1-e^{-\jf_w^{w+t}\l(u)du}$. Dividing both sides of (\ref{5}) and then considering $h\downarrow0$, due to (\ref{I1})-(\ref{I3}) and $[V-\vf](s,x,w)=0$, we obtain
\[
\zb[\!-cV\!\!+\!\vf_s\!+\!p\vf_x\!+\!\vf_w\yb]\!(s,x,w)\!+\!l_0(1\!-\!\vf_x)\!(s,x,w)
\!\!+\!\!\l(w)\!\zb[\jf_0^x\!\!V(s,x\!-\!u,0)dG(u)\!\!-\!\!V(s,x,w)\yb]\!\!\le\!\!0.
\]
As $l_0$ is arbitrary, we obtain
\[
\max\left\{{1-\varphi_x,\hx{L}(V,\varphi)}\right\}\le 0.
\]
At this point, we have shown that $V$ is a viscosity supersolution on $\hx{D}$.

 $Subsolution$ \q Now, we prove that $V$ is the viscosity subsolution of the HJB equation on $\hx{D}^*$. If we assume the contrary, then there exists a point $(s,x,w)\in\hx{D}^*$
 and $\psi^0\in \mathbb{C}^{1,1,1}({D})$ such that $0=[V-\psi^0](s,x,w)=\max_{(t,y,v)\in{\hx{D}^*}}[V-\psi^0](t,y,v)$, but
\[
\max\left\{1-\psi_x^0,\hx{L}(\psi^0)\right\}(s,x,w)=-2\eta<0,
\]
where $\eta>0$ is a constant and $\hx{L}$ is the first-order integro-differential operator. Now, we show that there exists a function $\psi\in{\mathbb{C}}^{1,1,1}(\dqy^*)$ and constants $\varepsilon>0$ and $\rho>0$ such that
\bea
\begin{split}\label{sub1}
 &\hx{L}[\psi](t,y,v)\le -\varepsilon c,  (t,y,v)\in \overline {B_\rho(s,x,w)\cap\hx{D}^*}\setminus \{t=T\};\\
    &1-\psi_y(t,y,v)\le -\varepsilon c,  (t,y,v)\in \overline {B_\rho(s,x,w)\cap\hx{D}^*}\setminus \{t=T\};\\
    &V(t,y,v)\le \psi (t,y,v)-\varepsilon, (t,y,v)\in\partial B_\rho(s,x,w)\cap \hx{D}^*,
\end{split}
\eea
where $B_\rho(R)$ denotes the open sphere centered at $R$ with radius $\rho$.

\no
$Case\;1$ \q If $x>0$,  define
\[
\psi(t,y,v)\dy\psi^0(t,y,v)+\frac{\eta\left[(t-s)^2+(y-x)^2+(v-w)^2\right]^2}{\lambda(w)(x^2+w^2)^2}.
\]
Clearly, $\p\in \mathbb{C}^{1,1,1}(D)$, $\psi(s,x,w)=\psi^0(s,x,w)=V(s,x,w)$ and $\psi(t,y,v)>V(t,y,v)$ for all $(t,y,v)\neq (s,x,w)$. Furthermore, $(\psi_t^0,\psi_y^0,\psi_v^0)(s,x,w)=(\psi_y,\psi_y,\psi_v)(s,x,w)$ and
$
\lambda(w)\int_0^x\psi(s,x-u,0)dG(u)\le \lambda(w)\int_0^x\psi^0(s,x-u,0)dG(u)+\eta.
$
Thus, we have
\[\hx{L}\!\zb[\psi\yb](s,x,w)\!\!\le\!\!\hx{L}\!\zb[\p^0\yb]\!(s,x,w)\!\!+\!\!\eta\!=\!-\!\eta\!\!<\!0, \quad
[1\!-\!\psi_x](s,x,w)\!\!=\!\![1\!-\!\psi_x^0](s,x,w)\!=\!-2\eta\!<\!0.\]
This leads to
\[\max\{1-\psi_x,\hx{L}[\psi]\}(s,x,w)\le-\eta<0.\]
By continuity of $\psi_x$ and $\hx{L}[\psi]$, we know that there exists a constant $\rho>0$, such that
\[\max\{1-\psi_x,\hx{L}[\psi]\}(t,y,v)\le-\frac{\eta}{2}<0,\quad \mbox{for}\; (t,y,v)\in\overline{B_\rho(s,x,w)\cap \hx{D}^*}\setminus \{t=T\}.\]
Note that for all $(t,y,v)\in\partial B_\rho(s,x,w)\cap \hx{D}^*$, the following holds:
\[V(t,y,v)\le \psi(t,y,v)-\frac{\eta \rho^4}{\lambda(w)(x^2+w^2)^2}.\]
If we choose $\varepsilon=\min\{\frac{\eta \rho^4}{\lambda(w)(x^2+w^2)^2},\frac{\eta}{2c}\}$, we obtain (\ref{sub1}).

$Case \;2$ \q If $x=0$, in this case,
\[ \max\!\left\{1\!-\!\psi_y^0,\hx{L}[\psi^0]\right\}\!(s,0,w)\!=\!\max\!\left\{1\!-\!\psi_y^0,-c\psi^0\!+\!\psi_s^0\!+\!p\psi_y^0\!+\!\psi_w^0\!-\!\lambda(w)\psi^0\right\}\!(s,0,w)\!\le\! -2\eta\!<\!0.\]
Define $\psi(t,y,v)\dy\psi^0(t,y,v)+\eta[(t-s)^2+y^2+(v-w)^2].$
If we denote $\varepsilon=\min\{\frac{\eta}{2c},\eta\rho^2\}$ and perform a calculation similar to that of Case 1, we can show that (\ref{sub1}) still holds. In the following, we will argue that (\ref{sub1}) leads to a contradiction.

For any strategy $\pi\in{U}_{ad}^{x,w}[s,T]$, denote $R_t^{s,x,w}=(t,X_t^{s,x,w},W_t^{s,x,w})$. Define $\tau_{\rho}\triangleq\inf\{t>s:R_t\notin \overline{{B_\rho (s,x,w)\cap \hx{D}^*}}\}$ and $\theta\triangleq\tau_\rho \wedge T_1^{s,w};$ additionally, we denote $R_t=R_t^{s,x,w}$ for simplicity.
Applying  It$\hat{\mbox{o}}$'s formula to $e^{-c(t-s)}\psi(R_t)$, we have
\begin{eqnarray}\nonumber
\psi (s,x,w) &=& \mathbb{E}_{sxw}\left[e^{-c(\theta-s)}\psi (\theta,X_\theta, W_\theta)\right]+ \mathbb{E}_{sxw}\left[\int_{s}^{\th}e^{-c(t-s)}\psi_xdL_t^c\right]\\\nonumber & & +\mathbb{E}_{sxw}\left[\int_{s}^\th e^{-c(t-s)}(c\psi -\psi_t-\psi_xp-\psi_w)(t,X_{t-},W_{t-})dt\right]\\\nonumber & &+\mathbb{E}_{sxw}\left[\int_{s}^{\th}\int_0^{X_{t-}}(\psi(t-,X_{t-},W_{t-})-\psi(t,X_{t-}-u,0))dG(u)\lambda(W_t)dt\right]\\\nonumber & &+\mathbb{E}_{sxw}\left[\sum_{s\le t<\theta}e^{-c(t-s)}(\psi(R_t)-\psi(R_{t+}))\right],
\end{eqnarray}
where $L_t^c$ denotes the continuous part of process $L_t$.
For time $t\in [s,\theta)$, $-\hx{L}[\psi](t,X_t,W_t)\ge c\varepsilon, \psi_y(t,X_t,W_t)\ge 1+c\varepsilon.$  The above equation can be transformed as follows:
\begin{eqnarray}\label{a66}
\psi (s,x,w)\!\ge\! \mathbb{E}_{sxw}\!\left[e^{-c(\theta -s)}\psi (\theta, X_\theta,W_\theta)\right]\!\!+\!\!\mathbb{E}_{sxw}\!\left[\int_{s}^{\theta}e^{-c(t-s)}c\varepsilon dt\right]\!\!+\!\mathbb{E}_{sxw}\!\left[\int_{s}^\th e^{-c(t-s)}dL_t\right].
\end{eqnarray}
On $\{\tau_{\rho}< T_1^{s,w}, (\theta, X_\theta,W_\theta)\neq (\theta,X_{\theta+},W_{\theta})\}$,  meaning that at time $\theta$, the wealth process jumps out of $\overline{B_\rho (s,x,w)\cap \hx{D}^*}$ due to dividend payments. There exists a random variable $\nu \in[0,1]$ and a point $(\th,X_\nu,W_{\theta})\in \partial B_\rho (s,x,w)\cap \hx{D}^*$ such that
\[
X_\nu=X_\theta-\nu [L_{\theta+}-L_\theta].
\]
From (\ref{sub1}), we observe that
\begin{equation}\label{a60}
\begin{split}
\psi(\theta,X_\nu,W_\theta)&\ge V(\theta,X_\nu,W_\theta)+\e.\\
\psi(\theta,X_\theta,W_\theta)-\psi(\theta,X_\nu,W_\theta)&\ge (c\varepsilon+1)(X_\theta-X_\nu)=(c\varepsilon+1)\nu(L_{\theta+}-L_\theta).
\end{split}
\end{equation}
From proposition \ref{a38}-(1), we obtain
\begin{align}\label{a61}  V(\theta,X_\nu,W_\theta) \ge V(\th,X_{\th+},W_\th)+(1-\nu)[L_{\th+}-L_\th].
\end{align}
Combining (\ref{a60}) with  (\ref{a61}), we observe that  on $\zb\{\tau_{\rho}< T_1^{s,w}, (\theta, X_\theta,W_\theta)\neq (\theta,X_{\theta+},W_{\theta})\yb\}$,
\begin{eqnarray}\label{a62}
\psi(\theta,X_\theta,W_\theta)\ge  V(\theta,X_{\theta+},W_\theta)+\varepsilon+(L_{\theta+}-L_\theta).
\end{eqnarray}
On the other hand,
in the case of $\{\tau_\rho <T_1^{s,w},(\theta,X_\theta,W_\theta)=(\theta,X_{\theta+},W_{\theta})\}$, meaning that at time $\th,$ $(\theta,X_\theta,W_\theta)\in\partial B_\rho (s,x,w)\cap \hx{D}^*$;  we obtain
\begin{eqnarray}\label{a63}
\psi(\theta,X_\theta,W_\theta)\ge V(\theta,X_\theta,W_\theta)+\varepsilon=V(\theta,X_{\theta+},W_{\theta})+\varepsilon.
\end{eqnarray}
In the case of $\{\tau_\rho\ge T_1^{s,w},X_{T_1^{s,w}}<0\}$, we have
\begin{align}\label{a64}
\psi(\theta,X_\theta,W_\theta)=V(\theta,X_\theta,W_\theta) =0.
\end{align}
In the case of $\{\tau_\rho\ge T_1^{s,w},X_{T_1^{s,w}}\ge 0\}$, we have
 \begin{align}\label{a65}
 \psi(\theta,X_\theta,W_\theta)\ge V(\theta,X_\theta,W_\theta)\ge V(\theta+,X_{\theta+},W_{\theta+})+L_{\theta+}-L_{\theta}.
 \end{align}
 From (\ref{a62}), (\ref{a63}), (\ref{a64}), and (\ref{a65}), we obtain
 \begin{align}\bbh
 &\mathbb{E}_{sxw}\left[\!e^{-c(\theta-s)}\psi(\theta,X_\theta,W_\theta)+\!\int_{s}^{\th}e^{-c(t-s)}dL_t\right]\\\bbh
 \ge&\mathbb{E}_{sxw}\left[e^{-c(\theta-s)}V(\theta,X_{\theta+},W_{\theta})\!+\!\int_{s}^{\theta+}\!\! e^{-c(t-s)}dL_t:\tau_\rho\!\ge\! T_1^{s,w}\right]\\\label{22}
 &+\mathbb{E}_{sxw}\left[e^{-c(\theta-s)}V(\theta,X_{\theta+},W_{\theta})\!+\!\!\int_{s}^{\th+} e^{-c(t-s)}dL_t\!+\!\varepsilon e^{-c(\theta-s)}:\tau_\rho< T_1^{s,w}\right].
 \end{align}
Combining (\ref{a66})  with (\ref{22}), we observe that
\begin{align}\bbh
\psi(s,x,w)&\ge\! \mathbb{E}_{sxw}\!\left[e^{-c(\theta-s)}V(\theta,X_{\theta+},W_{\theta})\!+\!\!\!\int_{s}^{\th+}\!\!e^{-c(t-s)}dL_t\right]\!\!+\!\varepsilon\!-\!\varepsilon \mathbb{E}_{sxw}\left[e^{-c(\theta-s)}:\tau_\rho \!\ge\! T_1^{s,w}\right]\\\label{a67}
&\ge \mathbb{E}_{sxw}\!\left[e^{-c(\theta-s)}V(\theta,X_{\theta+},W_{\theta})\!+\!\!\int_{s}^{\th+}e^{-c(t-s)}dL_t\right]\!+\!\varepsilon \mathbb{E}_{sxw}\!\left[1\!-\!e^{-c(T_1^{s,w}-s)}\right].
\end{align}
For $h>0$, $\theta+h$ is a stopping time.
The dynamic programming principle yields that
\[
V(s,x,w)=\sup_{\pi\in{U}_{ad}^{x,w}[s,T]}\mathbb{E}_{sxw}\left[\int_{s}^{\th+h} e^{-c(t-s)}dL_t+e^{-c((\theta+h)\wedge \tau^\pi-s)}V(R_{(\theta+h)\wedge \tau^\pi})\right].
\]
Let $h\downarrow 0$; by continuity of $V$, we obtain
\begin{eqnarray}\label{23}
V(s,x,w)=\sup_{\pi\in{U}_{ad}^{x,w}[s,T]}\mathbb{E}_{sxw}\left[\int_{s}^{\th+}e^{-c(t-s)}dL_t\right]+e^{-c(\theta-s)}
V(\theta,X_{\theta+},W_{\theta}).
\end{eqnarray}
As $\mathbb{P}(T_1^{s,w}>s)=1$ and  $V(s,x,w)=\p(s,x,w)$, we observe that (\ref{23}) contradicts  (\ref{a67}).   $\hfill\Box$
\section{The Candidate of The  Value Function}
In this section, we provide a candidate of the value function.  Before we state theorems, we introduce several notations.
\begin{defn}
If $u: D\yjt \mathbb{R}^3$, define
\[u^*(s,x,w)=\lim_{r\downarrow 0}\sup\left\{u(t,y,v):(t,y,v)\in D \;\bh{and} \;\sqrt{|t-s|^2+|y-x|^2+|v-w|^2}\le r \right\}.\]
We call $u^*$ the upper semicontinuous envelope of $u$. $u^*$ is the smallest upper semi-continuous function satisfying $u\le u^*$.
\end{defn}
 Denote $d_{\dqy}(s,x,w)$ the distance between $(s,x,w)$ and the boundary of $\hx{D}$, which means
\begin{equation}\label{juli1}
d_{\dqy}(s,x,w)=(T-s)\qx w\qx \fs{\sqrt{2}}{2}(s-w)\qx x.
\end{equation}
 Recall the definition
$\hx{D}=\bh{int} D=\{(s,x,w)\in D| 0<s<T,0<x, 0<w<s\}$.
We observe that $d_\dqy(s,x,w)\le\fs{T}{2+\sqrt{2}}$. In fact, if we consider $w=T-s=\fs{\sqrt{2}}{2}(s-w)$, we obtain $w=T-s=\fs{\sqrt{2}}{2}(s-w)=\fs{T}{2+\sqrt{2}}.$ If
$x\qx w\qx T-s\qx \fs{\sqrt{2}}{2}(s-w)>\fs{T}{2+\sqrt{2}}$, then by summing $w,s-w$ and $T-s$ we obtain $T>T$, which is a contradiction.

 From now on,  denote $M_1\dy(\fs{\gh{2}}{2}+1+2p)T+\fs{T}{2+\gh{2}}$ and $M_2\dy-[\fs{\gh{2}}{2}+1+(c+\L)\fs{2T}{2+\gh{2}}]T$ for simplicity.
In what follows, we construct a viscosity supersolution of (\ref{1}) on $\hx{D}$.
\begin{thm} Define
\bea\label{20}
\overline{V}(s,x,w)=x+{d_{\dqy}(s,x,w)}+N_1(T-s),
\eea
where the constant $N_1=\fs{\gh{2}}{2}+1+2p$, and $d_{\dqy}$ is defined in (\ref{juli1}),
Then, $\ol{V}$
is a viscosity supersolution of (\ref{1}) on $\hx{D}$ and for all $(s,x,w)\in D$, $x+M_2\le\ol{V}(s,x,w)\le x+M_1$.
\end{thm}
$\mathbf{Proof}$\q Note that for all $(s,x,w)\in \hx{D}$, $d_{\hx{D}}(s,x,w) \le\fs{T}{2+\sqrt{2}}$. Thus, it is obvious that for all $(s,x,w)\in D$, $ x+M_2\le\ol{V}(s,x,w)\le x+M_1$.
Denote $R\dy (s,x,w)$ for simplicity. Consider a function $\f\in \mathbb{C}^{1,1,1}(\hx{D})$ such that $\ol{V}-\f$ attains its minimum at $R$, which means $x+{d_{\dqy}(s,x,w)}+N_1(T-s)-\f$ attains its minimum at $R$. Now, we discuss various cases.

\no
Case 1\q
If $\fs{\gh{2}}{2}(s-w)<w\qx x\qx (T-s)$, we observe that $d_{\hx{D}}(s,x,w)=\fs{\gh{2}}{2}(s-w)$ near point $R$; this leads to $(\f_x,\f_s,\f_w)(R)=(1,\fs{\gh{2}}{2}-N_1,-\fs{\gh{2}}{2})$.  Thus, we obtain
  \begin{align}\bbh
    \hx{L}[\ol{V}, \f](R)=&\zb[-(c+\l(w))\ol{V}+\f_s+p\f_x+\f_w\yb](R)
 +\l(w)\jf_0^{x}\ol{V}(s,x-u,0)dG(u)\\ \bbh
 \le& -c\zb[x\!+\!d_{\dqy}(R)\!+\!N_1(T\!-\!s)\yb]\!-\!\l(w)\zb[x\!+\!d_{\dqy}(R)\!+\!N_1(T\!-\!s)\yb]\!+p\!\\\bbh
 &-N_1\!+\!\l(w)\jf_0^{x}(x-u)dG(u)\!+\!\l(w)N_1(T\!-\!s)G(x)
 \!\le\!-\!N_1\!\!+\!p<0.
\end{align}
 It is obvious that $\max\zb\{1-\f_x,\hx{L}[\ol{V},\f]\yb\}(R)\le 0.$

\no
Case 2\q
If $w<\fs{\gh{2}}{2}(s-w)\qx x\qx (T-s)$, then $\f_w(R)=1, \f_s(R)=-N_1, \f_x(R)=1.$ Similarly to Case 1, we can deduce that $\max\{1-\f_x,\hx{L}[\ol{V},\f]\}(R)\le -N_1+p+1<0 $. Similarly, we can verify the cases of $x<\fs{\gh{2}}{2}(s-w)\qx w\qx (T-s)$ and $T-s<\fs{\gh{2}}{2}(s-w)\qx w\qx x$.

\no
Case 3\q
If  $\fs{\gh{2}}{2}(s-w)=w< (T-s)\qx x$, then $1\le \f_w(R)\le -\fs{\gh{2}}{2}$. Thus, we conclude that $\f$ does not exist. Similarly, we can verify that $\ol{V}$ is a viscosity supersolution of \ref{1} in other cases of $\hx{D}$. $\hfill\Box$

In what follows, we provide a subsolution of (\ref{1}) on $\hx{D}^*$.
\begin{thm}
Define \begin{equation}\label{a12}
\underline{V}(s,x,w)=x+{d_{\hx{D}}(s,x,w)}-N_2(T-s),
\end{equation}
where the constant $N_2= \fs{\gh{2}}{2}+1+(c+\L)\fs{2T}{2+\gh{2}}$,
and $d_\dqy(s,x,w)$ is defined in (\ref{juli1}). Then, $\ul{V}$ is a viscosity subsolution of (\ref{1}) on $\dqy^*$, and for all $(s,x,w)\in D$, $x+M_2\le\ul{V}(s,x,w)\le x+M_1$.
\end{thm}
$\mathbf{Proof}$ \q It is obvious that for all $(s,x,w)\in D$, $x+M_2\le\ul{V}(s,x,w)\le x+M_1$. Now, we show that $\ul{V}$ is a viscosity subsolution of (\ref{1}) on $\hx{D}^*$. For a fixed point $R\dy(s,x,w)\in \dqy^*$, consider the function $\vf\in \mathbb{C}^{1,1,1}(D)$ such that $\ul{V}-\vf$ reaches its maximum at $R.$

\no Case 1\q $x>\fs{T}{2+\sqrt{2}}$.
First, we note  for all  $x>\fs{T}{2+\sqrt{2}}$, $x>\max\{T-s,w,\fs{\gh{2}}{2}(s-w)\}$, which means that for all $x>\fs{T}{2+\sqrt{2}}$, $\vf_x(R)=1$. Obviously, this  leads to $\max\zb\{1-\f_x,\hx{L}[\ol{V},\vf]\yb\}(R)\ge 0.$  Thus, we have shown  that $\ul{V}$ is a viscosity subsolution for $R\in [0,T]\times(\fs{T}{2+\sqrt{2}},\wq)\times [0,s]$.

\no Case 2 \q
If $R\in int\hx{D}^*$ and $x< (T-s)\qx w\qx\fs{\sqrt{2}}{2}(s-w)$, it is obvious that $x$ is less than $\fs{T}{2+\sqrt{2}}$. By several simple calculations, we obtain $\vf_x(R)=2,\vf_s(R)=N_2,\vf_w(R)=0.$ Then,
\begin{align}
\bbh
\hx{L}[\ul{V},\vf](R)=&\zb[-(c+\l(w))\ul{V}+\vf_s+p\vf_x+\vf_w\yb](R)
 +\l(w)\jf_0^{x}\ul{V}(s_0,x-u,0)dG(u)\\\bbh
\ge &-(c+\l(w))\zb[\fs{T}{2+\sqrt{2}}+\fs{T}{2+\gh{2}}\yb]+N_2+2p>0.
\end{align}
Thus, $\ul{V}$ is a viscosity subsolution of (\ref{1}) at $R$.  Now, we prove a special case of the boundary.

\no Case 3\q
 In the case of  $R=(0,0,0)\in\pw\hx{D}^*$, we observe that $\vf_s(R_0)\ge N_2, [\vf_s+\vf_w](R_0)\ge N_2,\vf_x(R_0) \ge1.$ Thus, it is easy to verify that
$\hx{L}\zb[\ul{V},\vf\yb](R)\ge N_2+p>0.$
Now, we observe that $\ul{V}$ is a viscosity subsolution of (\ref{1}) at $(0,0,0)$. We can prove other cases similarly.$\hfill\Box$
\begin{lem}  \label{a3} Let $\ol{V}$ be the viscosity supersolution constructed in  (\ref{20}), and $\underline{V}$ be the  viscosity subsolution constructed in (\ref{a12}). Define
\bea\label{supre}
\omega(s,x,w)=\sup_{u\in\hx{G}}u(s,x,w),
\eea
where $\hx{G}\dy$
\{$u|u$ is a viscosity subsolution of (\ref{1}) on $\hx{D^*}$ and $\ul{V}\le u\le\ol{V}$ in $D$\}.
Then,   $x+M_2\le\omega^*(s,x,w)\le M_1+x$, and $\omega^*$ is a viscosity subsolution of (\ref{1}) on $\hx{D^*}$.
\end{lem}
$\mathbf{Proof}$ \q   It is not difficult to verify that for all $(s,x,w)\in D$, $x+M_2\le\omega^*(s,x,w)\le M_1+x$.  In what follows, we borrow several arguments from \cite{mou2017perron}.   Suppose that there exists a $\vf\in\mathbb{C}^{1,1,1}(\hx{D}^*)$  such that $\omega^*-\vf$ attains its maximum (equal $0$) at $(s_0,x_0,w_0)\in \hx{D}^*$ over $D$. There exists $\eta_2>\eta_1>0$ and a  sequence of $\mathbb{C}^{1,1,1}(D)$ functions $\{\vf_m\}_{m}$ such that $\vf_m=\vf$  in $B_{\eta_1}(s_0,x_0,w_0)$, $\vf\le \vf_m$ in $\hx{D}^*$,
\[\sup_{(s,x,w)\in B_{\eta_2}(s_0,x_0,w_0)}\{\omega^*(s,x,w)-\vf_m(s,x,w)\}\le-\fs{1}{m}\]
and $\vf_m\yjt\vf$ pointwise. Thus, for any positive integer $m$, $\omega^*-\vf_m$ has a strict maximum of 0 at $(s_0,x_0,w_0)\in \hx{D}^*$ over $D$. Therefore,
\[
\sup_{(s,x,w)\in B_{\eta_1}^c(s_0,x_0,w_0)\cap D}\{\omega^*(s,x,w)-\vf_m(s,x,w)\}= \e_m<0.
\]
By the definition of $\omega^*$, we have, for any $u\in\hx{G}$,
\[\sup_{(s,x,w)\in B_{\eta_1}^c(s_0,x_0,w_0)\cap {D}}\{u(s,x,w)-\vf_m(s,x,w)\}\le \e_m<0.\]
Again, by the definition of $\omega^*$, we note that for any $\e<0$ satisfying  $\e>\e_m$, there exists $u_\e\in\hx{G}$ and $(\bar{s}_\e,\bar{x}_\e,\bar{w}_\e)\in B_{\eta_1}(s_0,x_0,w_0)\cap D$  such that
\begin{align}\bbh u_\e(\bar{s}_\e,\bar{x}_\e,\bar{w}_\e)-\vf(\bar{s}_\e,\bar{x}_\e,\bar{w}_\e)>\e.\end{align}
Since $u_\e\in USC({D})$  and $\vf_m\in \mathbb{C}^{1,1,1}(D)$, there exists $(s_\e,x_\e,w_\e)\in B_{\eta_1}(s_0,x_0,w_0)\cap D$ such that \begin{align}\label{a9}
u_{\e}(s_\e,x_\e,w_\e)-\vf_m(s_\e,x_\e,w_\e)=\sup_{(s,x,w)\in D}\{u_\e(s,x,w)-\vf_m(s,x,w)\}.
\end{align}
Since $\vf_m=\vf$ on $B_{\eta_1}(s_0,x_0,w_0)$, we deduce that
\begin{align}\bbh
\sup_{(s,x,w)\in{D}}\{u_\e(s,x,w)-\vf_m(s,x,w)\}&\ge\sup_{(s,x,w)\in B_{\eta_1}(s_0,x_0,w_0)\cap D}\{u_\e(s,x,w)-\vf(s,x,w)\}\\\label{a10}
&\ge u_{\e}(\bar{s}_\e,\bar{x}_\e,\bar{w}_\e)-\vf(\bar{s}_\e,\bar{x}_\e,\bar{w}_\e)>\e.
\end{align}
Combining (\ref{a9}) and (\ref{a10}), we obtain
\[
u_{\e}(s_\e,x_\e,w_\e)-\vf_m(s_\e,x_\e,w_\e)>\e.
\]
Since $\omega^*-\vf_m$ attains its strict maximum of $0$ over $\hx{D}^*$ and $\f\le \omega^*$ for any $\f\in\hx{G}$, then $u_{\e}(s_\e,x_\e,w_\e)\yjt \omega^*(s_0,x_0,w_0)$  and $(s_\e,x_\e,w_\e)\yjt (s_0,x_0,w_0)$ as $\e\yjt 0^-$. Denote $C_{m,\e}=u_{\e}(s_\e,x_\e,w_\e)-\vf_m(s_\e,x_\e,w_\e)$.  As $u_\e$ is a viscosity subsolution of the HJB equation (\ref{1}) and $u_\e-(\vf_{m}+C_{m,\e})$ attains its strict maximum (equal $0$) at point $(s_\e,x_\e,w_\e)$ in $\hx{D}^*$, we have
\begin{align}\bbh
\max\bigg\{&[-c(\vf_m\!+\!C_{m,\e})\!+\!\{\vf_m\}_s\!+\!p\{\vf_m\}_x\!+\!\{\vf_m\}_w](s_\e,x_\e,w_\e)\!-\!\l(w_\e)\zb[\vf_m(s_\e,x_\e,w_\e)\!+\!C_{m,\e}\yb]\\\bbh
&+\!\l(w_\e)\jf_0^{x_\e}\vf_m(s_\e,x_\e\!-\!u,0)dG(u),\qquad 1\!-\!\{\vf_{m}+C_{m,\e}\}_x\bigg\}\ge 0.
\end{align}
Letting $\e\yjt0^-$ and $m\yjt+\wq$, we have
\begin{align}\bbh
\max\bigg\{&[-c\vf\!+\!\vf_s\!+p\vf_x\!+\!\vf_w](s_0,x_0,w_0)\!-\!\l(w_0)\vf(s_0,x_0,w_0)\\\bbh
&+\l(w_0)\jf_0^{x_0}\vf(s_0,x_0-u,0)dG(u),\quad 1-\vf_x\bigg\}\ge 0.
\end{align}
Here, we used the fact that $C_{m,\e}\yjt 0$, as $\e\yjt 0^-$ and $m\yjt\wq$. This completes the proof. $\hfill\Box$

At this point, we are ready to provide the representation of the value function $V$.
\begin{thm}
The value function $V=\omega^*$, where $\o$ is defined in (\ref{supre}).
\end{thm}
$\mathbf{Proof}$\q From Lemma \ref{a3}, we see that $\omega^*$ is a viscosity subsolution of (\ref{1}). Note that the value function $V$ is a viscosity subsolution of (\ref{1}) that satisfies $\ul{V}\le V\le\ol{V}$.  By the definition of $\omega^*$, we know that $\omega^*\ge V$. Now, we only need to show that $V\ge \omega^*$.
Choose a sufficiently large $\yw{K}$ such that  $c\yw{K}\ge 1.$  Consider the function
\[
V^{\theta}=\th V+(\th-1)\yw{K},
\]
where $\th>1$.
It is easy to verify that $V^{\th}$ is also a  viscosity supersolution of (\ref{1}). Indeed, consider any continuously differentiable function $\vf$ such that $V^{\th}-\vf$ attains its minimum at $R_0$, which means $V-\zb[\fs{\vf}{\th}-\fs{\th-1}{\th}\yw{K}\yb]$ attains its minimum at $R_0\in \hx{D}$. Then,
\[
\max\zb\{1-\zb(\fs{\vf}{\th}-\fs{\th-1}{\th}\yw{K}\yb)_x,\hx{L}\zb[V,\fs{\vf}{\th}-\fs{\th-1}{\th}\yw{K}\yb]\yb\}\le0.
\]
Thus,
\begin{align}\label{4}
\max\zb\{1-\vf_x,\hx{L}\zb[V^\th,\vf\yb]\yb\}\le-(\th-1).
\end{align}
This shows that $V^{\th}$ is also a  viscosity supersolution of (\ref{1}). Instead of comparing $\omega^*$ and $V$, we will compare $\omega^*$ and $V^\th$. If we can show that $\omega^*\le V^\th$, then by simply sending $\th\yjt 1+$, we obtain the desired comparison result $\omega^*\le V$ in $D.$ Observe that for all $(s,x,w)\in D$, $x+M_2\le V(s,x,w)\le x+M_1$ and $x+M_2\le \omega^*(s,x,w)\le x+M_1$; thus, we have
\begin{align}\label{13}
[\omega^*\!-\!V^\th](s,x,w)\le x+M_1\!-\!\zb[\th(x\!+\!M_2)\!+\!(\th\!-\!1)\yw{K}\yb]
\!=\!(1-\th)x\!+\!(M_1-\th M_2-(\th-1)\yw{K}).
\end{align}
In view of (\ref{13}), we can choose $b\dy\fs{M_1-\th M_2-(\th-1)\yw{K}}{\th-1}$ such  that for all $x\ge b$, we have $\omega^*\le V^\th$. Although $D$ is bounded, we can then restrict our attention to the bounded domain
\[
\hx{D}_b\dy\zb\{(t,y,v)|0\le t<T,0\le y<b,0\le v\le t\yb\}
\]
and prove that $\omega^*\le V^\th$ on $\hx{D}_b$. Now, we assume on the contrary that
\[
M_b\dy\max_{\overline{\hx{D}_b}}\zb[\omega^*-V^{\th}\yb]=(\omega^*-V^\th)(\bar{s},\bar{x},\bar{w})>0
\]
for some $\bar{R}\dy(\bar{s},\bar{x},\bar{w})\in \ol{\hx{D}_b}$. Observe that we have only two cases $\bar{R}\in \bh{int}\hx{D}_b\dy\hx{D}_b^0$ and $\bar{R}\in \hx{D}_b^1$ to consider, where
\[
\hx{D}_b^1\dy\pw\hx{D}_b\qd[\{t=T\}\cup\{y=b\}]
\]
is the state constraint boundary restricted by $b.$

Case \uppercase\expandafter{\romannumeral1}\q Consider $\bar{R}\in \hx{D}_b^1$. The construction presented below is a suitable  adaption of the construction of \cite{benth2001optimal}. Denote $R\dy(t,y,v)$ for simplicity.  As $\hx{D}_b$ is  piecewise linear, there exist constants  $h_0,k>0$ and a uniformly  continuous map
$\eta:\ol{\hx{D}_b}\mapsto \mathbb{R}^3$ satisfying
\bea\label{16}
B(R+h\eta(R), hk)\ysy\hx{D}_b^0 \;\bh{for}\;\bh{all}\; R\in\ol{\hx{D}_b} \;\bh{and}\;{h}\in (0,h_0],
\eea
where $B(z,\r)$ denotes the sphere with radius $\r$ and center $z$. We can write   $\eta({R})=(\eta_1({R}),\eta_2({R}),\eta_3({R}))$.  For any $\k>1$ and $0<\e<1$, define the function $\Phi(s,x,w,t,y,v)$ on $\ol{\hx{D}_b}\times \ol{\hx{D}_b}$  by
\begin{align}\bbh
\Phi(s,x,w,t,y,v)=&\omega^*(s,x,w)-V^\th(t,y,v)-(\k (s-t)+\e\eta_1(\bar{R}))^2-(\k(x-y)+\e\eta_2(\bar{R}))^2\\\bbh
&-(\k(w-v)+\e\eta_3(\bar{R}))^2-\e(s-\bar{s})^2-\e(x-\bar{x})^2-\e(w-\bar{w})^2.
\end{align}
Let\[
M_\k=\max_{\ol{\hx{D}_b}\times\ol{\hx{D}_b}}\Phi(s,x,w,t,y,v).
\]
We then have  $M_\k\ge \omega^*(\bar{s},\bar{x},\bar{w})-V^\th(\bar{s},\bar{x},\bar{w})-\e^2|\eta(\bar{R})|^2>0$ for any $\k>1$ and $\e<\e_0$, where $\e_0$ is a small fixed number. Let $(s_\k,x_\k,w_\k,t_\k,y_\k,v_\k)\in \ol{\hx{D}_b}\times\ol{\hx{D}_b}$ be a  maximizer of $\Phi$, i.e., $M_\k=\Phi(s_\k,x_\k,w_\k,t_\k,y_\k,v_\k)$. By (\ref{16}), we assume that $\k$ is so large that $\bar{R}+\fs{\e}{\k}\eta(\bar{R})\in \hx{D}_b^0$. From
\[
\Phi(s_\k,x_\k,w_\k,t_\k,y_\k,v_\k)\ge \Phi(\bar{s},\bar{x},\bar{w},\bar{s}+\fs{\e}{\k}\eta_1(\bar{R}),\bar{x}+\fs{\e}{\k}\eta_2(\bar{R}),\bar{w}+\fs{\e}{\k}\eta_3(\bar{R})),
\]
we obtain that
\begin{align}\bbh
&|\k(s_\k\!-\!t_\k)\!+\!\e\eta_1(\bar{R})|^2\!\!+\!\!|\k(x_\k\!-\!y_\k)\!+\!\e\eta_2(\bar{R})|^2\!\!+\!\!|\k(w_\k\!-\!v_\k)\!+\!\e\eta_3(\bar{R})|^2\\\bbh
&+\e(s_\k\!-\!\bar{s})^2+\e(x_\k\!-\!\bar{x})^2\!+\!\e(w_\k-\bar{w})^2\\\bbh
\le& \omega^*(s_\k,x_\k,w_\k)-V^\th(t_\k,y_\k,v_\k)-(\omega^*-V^\th)(\bar{s},\bar{x},\bar{w})-V^\th(\bar{s},\bar{x},\bar{w})\\\label{17}
&+V^\th(\bar{s}+\fs{\e}{\k}\eta_1(\bar{R}),\bar{x}+\fs{\e}{\k}\eta_2(\bar{R}),\bar{w}+\fs{\e}{\k}\eta_3(\bar{R})).
\end{align}
As $\omega^*$ and $-V^\th$ are bounded on $\bar{\hx{D}_b}$, it follows that $|\k(s_\k-t_\k)|,|\k(x_\k-y_\k)|,|\k(w_\k-v_\k)|$ are bounded uniformly in $\k$. Hence, we have
\[
s_\k\!-\!t_\k\!\yjt\! 0,x_\k\!-\!y_\k\!\yjt\! 0,w_\k\!-\!v_\k\!\yjt\! 0, \;\bh{as}\; \k\yjt\wq\; \bh{and}
\lim_{\k\yjt\wq}\!\zb[\omega^*(s_\k,x_\k,w_\k)\!-\!V^\th(t_\k,y_\k,v_\k)\yb]\le M_b.
\]
Sending $\k\yjt\wq$ in (\ref{17}) and using the upper semi-continuity of $\omega^*$, $-V^\th$ in $\ol{\hx{D}_b}$, we conclude that
\[|\k(s_\k-t_\k)+\e\eta_1(\bar{R})|\yjt0,|\k(x_\k-y_\k)+\e\eta_2(\bar{R})|\yjt0, |\k(w_\k-v_\k)+\e\eta_3(\bar{R})|\yjt0,\]
$(s_\k,x_\k,w_\k)\yjt \bar{R}$, $(t_\k,y_\k,v_\k)\yjt \bar{R}$ and $M_\k\yjt M_b.$ Using the uniform continuity of $\eta$, we have that
\[t_\k=s_\k+\fs{\e}{\k}\eta_1(\bar{R})+o\zb(\fs{1}{\k}\yb)=s_\k+\fs{\e}{\k}\eta_1(s_\k,x_\k,w_\k)+o\zb(\fs{1}{\k}\yb),\]
Similarly, we have that
\[y_\k=x_\k+\fs{\e}{\k}\eta_2(s_\k,x_\k,w_\k)+o\zb(\fs{1}{\k}\yb), \q v_\k=w_\k+\fs{\e}{\k}\eta_3(s_\k,x_\k,w_\k)+o\zb(\fs{1}{\k}\yb).\]
Thus, we use (\ref{16}) to obtain $(t_\k,y_\k,v_\k)\in\hx{D}_b^0$ for sufficiently large $\k$. Now, define
\begin{align}\bbh
\f(s,x,w)=&V^\th(t_\k,y_\k,v_\k)+|\k(s-t_\k)+\e\eta_1(\bar{R})|^2+|\k(x-y_\k)+\e\eta_2(\bar{R})|^2\\\bbh
&+|\k(w-v_\k)+\e\eta_3(\bar{R})|^2+\e(s-\bar{s})^2+\e(x-\bar{x})^2+\e(w-\bar{w})^2.\\\bbh
\vf(t,y,v)=&\omega^*(s_\k,x_\k,w_\k)-|\k(s_\k-t)+\e\eta_1(\bar{R})|^2-|\k(x_\k-y)+\e\eta_2(\bar{R})|^2\\\bbh
&-|\k(w_\k-v)+\e\eta_3(\bar{R})|^2-\e(s_\k-\bar{s})^2-\e(x_\k-\bar{x})^2-\e(w_\k-\bar{w})^2.
\end{align}
By a direct calculation, we observe that
\beq\bbh
\begin{split}
\f_s(s_\k,x_\k,w_\k)\!=&\!2\k[\k(s_\k\!-\!t_\k)\!+\!\e\eta_1(\bar{R})]\!+\!2\e(s_\k\!-\!\bar{s});\;\;\;\;\vf_t(t_\k,y_\k,v_\k)\!=\!2\k[\k(s_\k\!-\!t_\k)\!+\!\e\eta_1(\bar{R})];\\
\f_x(s_\k,x_\k,w_\k)\!=&\!2\k[\k(x_\k\!-\!y_\k)\!+\!\e\eta_2(\bar{R})]\!+\!2\e(x_\k\!-\!\bar{x});\;\;\vf_y(t_\k,y_\k,v_\k)\!=\!2\k[\k(x_\k\!\!-\!y_\k)\!+\!\e\eta_2(\bar{R})];\\
\f_w(s_\k,x_\k,w_\k)\!=&\!2\k[\k(w_\k\!-\!v_\k)\!\!+\!\e\eta_3(\bar{R})]\!+\!2\e(w_\k\!-\!\bar{w});\;\vf_v(t_\k,y_\k,v_\k)\!=\!2\k[\k(w_\k\!\!-\!v_\k)\!+\!\e\eta_3(\bar{R})].\\
\end{split}
\endeq
As $V^\th-\vf$ reaches its minimum at $(t_\k,y_\k,v_\k)\in \hx{D}_b^0$, and $\omega^*-\f$ attains its maximum at $(s_\k,x_\k,w_\k)\in \hx{D}_b$, combining this with (\ref{4}), we observe that
\begin{align}\bbh
\max\bigg\{&1\!-\!2\k\zb[\k(x_\k\!-\!y_\k)\!+\!\e\eta_2(\bar{R})\yb],\q -(c\!+\!\l(v_\k))V^\th(t_\k,y_\k,v_\k)\!+\!2\k[\k(s_\k-t_\k)\!+\!\e\eta_1(\bar{R})]\\\bbh
&\q\q\q\qq\q\qq\qq+2p\k\zb[\k(x_\k\!-\!y_\k)\!+\!\e\eta_2(\bar{R})\yb]\!+\!2\k\zb[\k(w_\k\!-\!v_\k)\!+\!\e\eta_3(\bar{R})\yb]\\\label{18}
&\qq\qq\qq\qq\qq+\!\l(v_\k)\!\!\jf_0^{y_\k}\!\!V^\th(t_\k, y_\k\!-\!u,0)dG(u) \bigg\}\le-(\th-1).\\\bbh
\max\bigg\{&1\!-\!2\k[\k(x_\k\!-\!y_\k)\!+\!\e\eta_2(\bar{R})]-2\e(x_\k-\bar{x}),\; -(c\!+\!\l(w_\k))\omega^*(s_\k,x_\k,w_\k)\\\bbh
&\q+2\k[\k(s_\k\!-\!t_\k)\!+\!\e\eta_1(\bar{R})]\!+\!2\e(s_\k\!-\!\bar{s})+p\zb[2\k[\k(x_\k\!-\!y_\k)+\e\eta_2(\bar{R})]+2\e(x_\k-\bar{x})\yb]\\\label{19}
&\q+2\k[\k(w_\k-v_\k)\!+\!\e\eta_3(\bar{R})]\!+\!2\e(w_\k-\bar{w})\!+\!\l(w_\k)\jf_0^{x_\k}\omega^*\!(s_\k,x_\k-u,0)dG(u)\!\bigg\}\!\ge\!0.
\end{align}
Combining (\ref{18}) with (\ref{19}), we send $\k\yjt\wq$, $\e\yjt0$ to obtain the desired contradiction
\[[c+\l(\bar{w})](\omega^*-V^\th)(\bar{s},\bar{x},\bar{w})<\l(\bar{w})(\omega^*-V^\th)(\bar{s},\bar{x},\bar{w}).\]

\no Case \uppercase\expandafter{\romannumeral2} \q Let us consider the case $\bar{R}\in\hx{D}_b^0$. For any $\k>1$, define the function $\P$ on $\ol{\hx{D}_b}\times \ol{\hx{D}_b}$ by
\[\P(s,x,w,t,y,v)=\omega^*(s,x,w)-V^{\th}(t,y,v)-\fs{\k}{2}(s-t)^2-\fs{\k}{2}(x-y)-\fs{\k}{2}(w-v)^2.\]
Let $M_\k=\max_{\ol{\hx{D}_b}\times\ol{\hx{D}_b}}\P(s,x,w,t,y,v)$. We have $M_\k\ge M_b>0$ for all $\k>1$. Let $(s_\k,x_\k,w_\k,t_\k,y_\k,v_\k) $ be a maximizer, so that $M_\k=\P(s_\k,x_\k,w_\k,t_\k,y_\k,v_\k)$.
As $\ol{\hx{D}_b}\times\ol{\hx{D}_b}$  is compact, we can find a subsequence that may be assumed to be $(s_\k,x_\k,w_\k,t_\k,y_\k,v_\k)$ itself, such that $(s_\k,x_\k,w_\k,t_\k,y_\k,v_\k)\yjt(\yj{s},\yj{x},\yj{w},\yj{t},\yj{y},\yj{v})$. From
$M_\k\ge \P(\bar{s},\bar{x},\bar{w},\bar{s},\bar{x},\bar{w}),$ we obtain
\[\omega^*(s_\k,x_\k,w_\k)\!-\!V^{\th}(t_\k,y_\k,v_\k)\!-\!\fs{\k}{2}(s_\k\!\!-\!t_\k)^2\!-\!\fs{\k}{2}(x_\k\!\!-\!y_\k)^2\!-\!\fs{\k}{2}(w_\k\!-\!v_\k)^2
\!\ge\!\omega^*(\bar{s},\bar{x},\bar{w})\!-\!V^\th(\bar{s},\bar{x},\bar{w}).\]
Thus, we observe that
\begin{align}\label{21}
\!\!\fs{\k}{2}(s_\k\!\!-\!\!t_\k)^2\!\!+\!\!\fs{\k}{2}(x_\k\!\!-\!\!y_\k)^2\!\!+\!\!\fs{\k}{2}(w_\k\!\!-\!\!v_\k)^2\!\!\le\! \! \omega^*(s_\k,x_\k,w_\k)\!\!-\!\!V^{\th}(t_\k,y_\k,v_\k)
\!\!-\!\!\omega^*(\bar{s},\bar{x},\bar{w})\!\!+\!\!V^\th(\bar{s},\bar{x},\bar{w}).
\end{align}
As an upper semi-continuous functions attains its maximum on any compact set, we obtain that $\fs{\k}{2}(s_\k-t_\k)^2+\fs{\k}{2}(x_\k-y_\k)^2+\fs{\k}{2}(w_\k-v_\k)^2$ is bounded uniformly in $\k$. Thus, we obtain $x_\k-y_\k\yjt0$, $s_\k-t_\k\yjt0$, $w_\k-v_\k\yjt 0$ as $\k\yjt\wq$, which means $\yj{s}=\yj{t},\yj{w}=\yj{v},\yj{x}=\yj{y}$. Sending $\k\yjt\wq$ in (\ref{21}), we have  that
\[
\lim_{\k\yjt\wq}\!\zb[\!\fs{\k}{2}\!(s_\k\!-\!t_\k)^2\!+\!\fs{\k}{2}(x_\k\!-\!y_\k)^2\!+\!\fs{\k}{2}(w_\k\!-\!v_\k)^2\!\yb] \!\!+\!\omega^*(\bar{s},\bar{x},\bar{w})\!-\!V^\th\!(\bar{s},\bar{x},\bar{w})\!
\le\! \omega^*(\yj{s},\yj{x},\yj{w})\!-\!V^\th(\yj{s},\yj{x},\yj{w}).
\]
By the definition of $(\bar{s},\bar{x},\bar{w})$, we obtain $\yj{s}=\bar{s}, \yj{x}=\bar{x},\yj{w}=\bar{w}.$
As $(\bar{s},\bar{x},\bar{w})\in \hx{D}_b^0$, we observe that for  sufficiently large $\k$, $(s_\k,x_\k.w_\k),(t_\k,y_\k,v_\k)\in \hx{D}_b^0$. Note that for any given $\k$,
$\P(s,x,w,t,y,v)$  attains its maximum at $(s_\k,x_\k,w_\k,t_\k,y_\k,v_\k)$. Define the functions
\begin{align}\bbh
\f(s,x,w)&=V^\th(t_\k,y_\k,v_\k)+\fs{\k}{2}(x-y_\k)^2+\fs{\k}{2}(s-t_\k)^2+\fs{\k}{2}(w-v_\k)^2,\\\bbh
\p(t,y,v)&=\omega^*(s_\k,x_\k,w_\k)-\fs{\k}{2}(s_\k-t)^2-\fs{\k}{2}(x_\k-y)^2-\fs{\k}{2}(w_\k-v)^2.
\end{align}
We observe that $\zb[\omega^*-\f\yb](s,x,w)$ attains its maximum at $(s_\k,x_\k,w_\k)\in\hx{D}_b^0$, and $\zb[V^\th-\p\yb](t,y,v)$ attains its minimum at $(t_\k,y_\k,v_\k)\in\hx{D}_b^0$; $\omega^*$ is a viscosity subsolution of (\ref{1}), and $V^\th$ is a viscosity supersolution of (\ref{1}); noting that $\max\{1-\p_x,\hx{L}[V^\th,\p]\}\le-(\th-1)$ because of (\ref{4}), we obtain
\begin{align}\bbh
\max\bigg\{1-\k(x_\k-y_\k),& -(c\!-\!\l(w_\k))\omega^*(s_\k,x_\k,w_\k)+p\k(x_\k-y_\k)+\k(s_\k-t_\k)+\k(w_\k-v_\k)\\\label{24}
&+\l(w_\k)\jf_0^{x_\k}\omega^*(s_\k,x_\k-u,0)dG(u)\bigg\}\ge 0,\\\bbh
\max\bigg\{1-\k(x_\k-y_\k),&-\zb[c+\l(v_\k)\yb]V^\th(t_\k,y_\k,v_\k)+\k(s_\k-t_\k)+\k(w_\k-v_\k)\\\label{11}
&+p\k(x_\k-y_\k)+\l(v_\k)\jf_0^{y_\k}V^\th(t_\k,y_\k-u,0)dG(u)\bigg\}\le-(\th-1).
\end{align}
From (\ref{11}), we observe that $\lim_{\k\yjt\wq}\zb[1-\k(x_\k-y_\k)\yb]<0.$ Thus, combining this with (\ref{24}) we obtain that for $\k$ large enough,
\bea\label{14}
-\![c\!\!+\!\!\l(\!w_\k\!)]\omega^*(s_\k,\!x_\k,\!w_\k)\!\!+\!\!\k(s_\k\!\!-\!t_\k\!)
\!\!+\!\!\k(w_\k\!\!-\!v_\k)\!\!+\!\!p\k(x_\k\!\!-\!\!y_\k)\!\!+\!\!\l(\!w_\k\!)\!\!\jf_0^{x_\k}\!\!\!\!\!\omega^*(s_\k,x_\k\!\!-\!u,\!0)dG\!(u)\!\ge\! 0
\eea
and
\bea\label{15}
\!-\!\zb[c\!+\!\!\l(\!v_\k\!)\yb]\!V^\th\!(t_\k,y_\k,v_\k\!)\!\!+\!\!\k\!(s_\k\!\!-\!\!t_\k)\!\!+\!\!\k(w_\k\!\!-\!\!v_\k)\!
+\!p\k\!(x_\k\!\!-\!\!y_\k\!)\!\!+\!\!\l(v_\k)\!\jf_0^{y_\k}\!\!\!V^\th(t_\k,y_\k\!-\!u,0)dG(u)\!\!<\!0.
\eea
Combining (\ref{14}) with (\ref{15}) and letting $\k\yjt\wq$, we obtain
\[-[c+\l(\bar{w})]M_b=-[c+\l(\bar{w})](\omega^*-V^\th)(\bar{s},\bar{x},\bar{w})> -\l(\bar{w})M_b,\]
which is a contradiction.
 We have shown that $\omega^*\le V$, which leads to $\omega^* = V$. $\hfill\Box$
\begin{rem}
In fact, the proof of $\omega^*\le V$ implies the comparison principle.  The comparison principle shows that for all supersolutions $\bar{u}$ and subsolutions $\underline{u}$, if $\bar{u}$ and $\ul{u}$ satisfy $x+M_2\le\bar{u}(s,x,w)\le x+M_1$ and $x+M_2\le\ul{u}(s,x,w)\le x+M_1$ for all $(s,x,w)\in D$, then $\ul{u}\le \bar{u}$.
\end{rem}
\section{Future work}
In our paper, there is a $``\max"$ operator in the HJB equation.  While aiming to analyze the optimal control problem, we know that exploring the representations of solutions of integro-PDEs becomes crucial.  Gong, Mou, and Swiech \cite{gong2017stochastic} explored the following integro-PDE:
\begin{equation}\label{mou}
\inf_{u\in U}\{\hx{A}^uW(t,x)+\G(t,x,u)\}=0 \;\;\bh{in} \;Q
\end{equation}
with the boundary condition $W(t,x)=\Phi(t,x)$, where $Q$ is a bounded domain, $\hx{A}$ is the generator of a drifted L$\acute{\bh{e}}$vy process with Brownian motion, $\G$ and $\Phi$ are given functions,
and $U$ is the value domain of control $u$. In  \cite{gong2017stochastic}, the authors state that the solution of (\ref{mou}) can be  approximated by a sequence of solutions of HJB integro-PDEs that are non-degenerate, have a finite control set, more regular coefficients and smooth terminal-boundary values on slightly enlarged domains. The authors showed that such slightly perturbed HJB equations have classical solutions.
Taking the limit of solutions can yield the solution of (\ref{mou}). As our HJB equation includes a ``max'' operator, it is unlikely to show that slightly perturbed HJB equations have continuously differentiable solutions.
Thus, we did not approximate a solution of the HJB equation in our paper.
 The closest study to our HJB equation exploring the regularity of solutions of obstacle integro-differential operators is that of Caffarelli,  Ros-Oton and Serra \cite{caffarelli2017obstacle}. In \cite{caffarelli2017obstacle}, the authors consider the obstacle problem in $\mathbb{R}^n$
\begin{equation}\label{obstacle1}
\begin{split}
\min\{-{L}u, u-\vf\}&=0 \q  \bh{in}\; \mathbb{R}^n\\
\lim_{|x|\yjt\wq}u(x)&=0,
\end{split}
 \end{equation}
 where $L$ is an infinitesimal generator of a L$\acute{\bh{e}}$vy process, and $\vf$ is a given bounded differentiable function on  $\mathbb{R}^n$ that we call an obstacle. The authors showed that a solution of (\ref{obstacle1}) belongs to $\mathbb{C}^{1,s}$ near all regular points, where $s\in (0,1)$. Unfortunately, our HJB equation $u_x$ is ``blocked" by $1$ instead of $u$ being ``blocked" by $1$.
  In the future work, we will focus on the structure of the HJB equation and try to find the optimal dividend strategy.

\end{CJK*}

\begin{thebibliography}{2}
    \bibitem{albrecher2005distribution} Albrecher, H., M$\acute{\bh{a}}$rmol, M., Claramunt, M. M., \emph{On the distribution of dividend payments in a Sparre Andersen model with generalized Erlang(n) interclaim times}. Insurance Mathematics and  Economics. 37, no. 2, 324-334, (2005).
        \bibitem{albrecher2006non}  Albrecher, H., Hartinger, J.,  \emph{On the non-optimality of horizontal barrier strategies in the Sparre Andersen model}. HERMES International Journal of Computer Mathematics and Its Applications, 7, 1-14,  (2006).
            \bibitem{albrecher2009optimality} Albrecher, H., Thonhauser, S., \emph{Optimality results for dividend problems in insurance}. RACSAM-Revista de la Real Academia de Ciencias Exactas, Fisicas y Naturales. Serie A. Matematicas, 103, no. 2, 295-320, (2009)
            \bibitem{andersen1957collective} Andersen, E. S., \emph{On the collective theory of risk in case of contagion between claims}. Bulletin of the Institute of Mathematics  and its Applications,  12, no. 2, 275-279, (1957).
    \bibitem{Asmussen1997Controlled} Asmussen, S., Taksar, M., \emph{Controlled diffusion models for optimal dividend pay-out}. Insurance: Mathematics and Economics, 20, no. 1, 1-15, (1997).
     \bibitem{awatif1991equqtions} Awatif, S., \emph{Equqtions d'Hamilton-Jacobi du premier ordre avec termes int{\'e}gro-diff{\'e}rentiels: Partie 1: Unicit{\'e} Des solutions de viscosit{\'e}}. Communications in partial differential equations, 16 (6-7):1057-1074, (1991).
      \bibitem{azcue2005optimal}  	Azcue, P., Muler, N., \emph{Optimal reinsurance and dividend distribution policies in the Cram$\acute{\bh{e}}$r-Lundberg model}. Mathematical Finance, 15, no. 2, 261-308, (2005).
       \bibitem{azcue2010optimal}  Azcue, P., Muler, N., \emph{Optimal investment policy and dividend payment strategy in an insurance company}. The Annals of Applied Probability, 20, no. 4, 1253-1302,  (2010).
        \bibitem{bai2017optimal}  Bai, L. H., Ma, J., Xing, X. J.,   \emph{Optimal dividend and investment problems under Sparre Andersen model}. The Annals of Applied Probability, 27, no. 6, 3588-3632,  (2017).
         \bibitem{belhaj2010optimal} Belhaj, M., \emph{Optimal dividend payments when cash reserves follow a jump-diffusion process.} Mathematical Finance, 20, no.2, 313-325, (2010).
         \bibitem{benth2001optimal} Benth, F. E., Karlsen, K. H.,  Reikvam, K., \emph{Optimal portfolio selection with consumption and nonlinear integro-differential equations with gradient constraint: a viscosity solution approach.} Finance and Stochastics, 5, no. 3, 275-303, (2001).
           \bibitem{caffarelli2017obstacle} Caffarelli, L. A., Ros-Oton, X., Serra, J., \emph{Obstacle problems for integro-differential operators: regularity of solutions and free boundaries}. Inventiones Mathematicae, 208, no. 3, 1155-1211, (2017).
            \bibitem{crandall1992user} Crandall, M. G.,  Ishii, H., Lions, P. L., \emph{User's guide to viscosity solutions of second order partial differential equations}. Bulletin of the American Mathematical Society, 27, no. 1, 1-67, (1992).
             \bibitem{crandall1983viscosity} Crandall, M. G., Lions, P. L.,  \emph{Viscosity solutions of Hamilton-Jacobi equations}. Transactions of the American Mathematical Society, 277, no. 1, 1-42, (1983).
              \bibitem{de1957impostazione} De Finetti, B., \emph{Su un'impostazione alternativa della teoria collettiva del rischio}. Transactions of the XVth international congress of Actuaries, 2, 433-443, New York, (1957).
\bibitem{fleming2006controlled} Fleming, W. H., Soner, H. M., \emph{Controlled Markov processes and viscosity solutions}. Second edition. Stochastic Modelling and Applied Probability, 25. Springer, New York, (2006).
              \bibitem{gerber2006optimal}Gerber, H. U., Shiu, S. W., \emph{On optimal dividend strategies in the compound Poisson model.} North American Actuarial Journal, 10, no. 2, 76-93, (2006).
 \bibitem{gong2017stochastic} Gong, R., Mou, C., Swiech, A., \emph{Stochastic Representations for Solutions to Nonlocal Bellman Equations.} arXiv preprint arXiv:1709.00193, (2017).
 \bibitem{li2004class} 	Li, S.M., Garrido, J., \emph{On a class of renewal risk models with a constant dividend barrier.} Insurance: Mathematics and Economics, 35, no. 3, 691-701, (2004).
     \bibitem{Lions1983} Lions, P. L., \emph{Optimal control of diffusion processes and Hamilton--Jacobi--Bellman equations part 2: viscosity solutions and uniqueness}.
     Communications in partial differential equations, 8, no. 11, 1229-1276, (1983).
 \bibitem{Lions1983Optimalone} Lions, P. L., \emph{Optimal control of diffusion processes and Hamilton-Jacobi-Bellman equations. I. The dynamic programming principle and applications.} Communications in Partial Differential Equations, 8, no. 10, 1101-1174, (1983).
 \bibitem{mou2017perron} Mou, C., \emph{Perron's method for nonlocal fully nonlinear equations.} Analysis $\&$ PDE, 10, no. 5, 1227-1254, (2017).
 \bibitem{Protter1990Stochastic} Protter, P., \emph{Stochastic Integration and Differential Equations.} Applications of Mathematics (New York), 21, Springer-Verlag, Berlin, (1990).
 \bibitem{rolski1998stochastic} Rolski, T., Schmidli, H., Schmidt, V., Teugels, J. L., \emph{Stochastic Processes for Insurance and Finance.} Vol. 505, John Wiley $\&$ Sons, (2009).
     \bibitem{scheer2011optimal} 	Scheer, N., Schmidli, H., \emph{Optimal dividend strategies in a Cram$\acute{e}$r--Lundberg model with capital injections and administration costs}. European Actuarial Journal, 1, no. 1, 57-92, (2011).
     \bibitem{seydel2009existence} Seydel, R. C., \emph{Existence and uniqueness of viscosity solutions for QVI associated with impulse control of jump-diffusions}. Stochastic Processes and their Applications, 119, no. 10, 3719-3748, (2009).
         \bibitem{soner1986optimal} Soner, H. M., \emph{Optimal control with state-space constraint. II}. SIAM Journal on Control and Optimization, 24, no. 6, 1110-1122, (1986).
             \bibitem{yong1999stochastic} Yong, J., Zhou, X. Y., \emph{Stochastic controls: Hamiltonian systems and HJB equations}. 43, Springer-Verlag, New York, (1999).
\end{thebibliography}
\end{document}